\documentclass[12pt]{amsart}
\theoremstyle{plain}
\newtheorem{theorem}{Theorem}
\newtheorem{lemma}{Lemma}
\newtheorem{corollary}[theorem]{Corollary}
\newtheorem{proposition}[theorem]{Proposition}
\theoremstyle{definition}

\newtheorem{remark}{Remark}

\def\Const{{\rm Const\ }}

\def\EXP{{\mathbb{E}}}
\def\PROB{{\mathbb{P}}}

\def\reals{{\mathbb{R}}}

\def\bdef{\begin{Def}}
\def\endef{\end{Def}}
\def\bthm{\begin{theorem}}
\def\ethm{\end{theorem}}
\def\blm{\begin{lemma}}
\def\elm{\end{lemma}}
\def\bcor{\begin{corollary}}
\def\ecor{\end{corollary}}
\def\bprop{\begin{proposition}}
\def\eprop{\end{proposition}}
\def\brm{\begin{remark}}
\def\erm{\end{remark}}
\def\beq{\begin{eqnarray}}
\def\eneq{\end{eqnarray}}
\def\beal{\begin{aligned}}
\def\enal{\end{aligned}}

\def\al{\alpha}

\def\bralpha{{\bar\alpha}}

\def\sg{\sigma}

\def\dt{\delta}

\def\cL{{\mathcal{L}}}
\def\cM{{\mathcal{M}}}

\def\cR{{\mathcal{R}}}
\def\cS{{\mathcal{S}}}
\def\cT{{\mathcal{T}}}
\def\cW{{\mathcal{W}}}

\def\R{\mathbb R}

\def\Z{\mathbb Z}
\def\N{\mathbb N}

\def\hX{{\hat X}}

\def\ha{{\hat a}}

\def\hrho{{\hat\rho}}
\def\th{\theta}
\def\~{\tilde}
\def\lb{\lambda}

\author{D. Dolgopyat, V. Kaloshin, L. Koralov}
\title{
Sample Path Properties of the Stochastic Flows.}
\address{\begin{tabular}{lll}
Dmitry Dolgopyat & Vadim Kaloshin & Leonid Koralov \\
Department of Mathematics & Department of &
Department of Mathematics \\
PennState University & Mathematics MIT & Princeton University \\
University Park PA 16802 & Cambridge MA 02139 & Princeton NJ 08544 \\
dolgop@math.psu.edu & kaloshin@math.mit.edu &
koralov@math.princeton.edu \\
www.math.psu.edu\slash dolgop\slash  & & \\
\end{tabular}}
\thanks{D.D. was partially supported by NSF and Sloan Foundation,
V. K. was partially supported by American Institute of Mathematics
Fellowship and Courant Institute, and L. K. was partially supported
by NSF postdoctoral fellowship.}
\begin{document}
\begin{abstract}
We consider a stochastic flow driven by a finite dimensional
Brownian motion. We show that almost every realization of such a
flow exhibits strong statistical properties such as the
exponential convergence of an initial measure to the equilibrium
state and the central limit theorem. The proof uses new estimates
of the mixing rates of the multipoint motion.
\end{abstract}
\maketitle
\markboth{Sample Path Properties of the Stochastic
 Flows}{D. Dolgopyat, V. Kaloshin, L. Koralov}

\section{Introduction.}

The subject of this paper is the study of the  long-time behavior
of a passive substance (or scalar) carried by a stochastic flow.
Motivation comes from applied problems in statistical turbulence
and  oceanography, Monin \& Yaglom  \cite{MY}, Yaglom \cite{Y},
Davis \cite{Da}, Isichenko \cite{I}, and Carmona \& Cerou
\cite{CC}. The questions we discuss here are also related to the
physical basis of the Kolmogorov model for turbulence, Molchanov
\cite{Mo}.

The physical mechanism of turbulence is still not completely
understood. It was suggested in \cite{RT} that the appearance of
turbulence could be similar to the appearance of chaotic behavior
in finite-dimensional deterministic systems. Compared to other
situations, the mechanism responsible for stochasticity in
deterministic dynamical systems with nonzero Lyapunov exponents is
relatively well understood. It is caused by a sensitive dependence
on initial conditions, that is, by exponential divergence of
nearby trajectories. It is believed that a similar mechanism can
be found in many other situations, but mathematical results are
scarce. Here we describe a setting where analysis similar to the
deterministic dynamical systems with nonzero Lyapunov exponents
can be used. In the paper we shall consider a  flow of
diffeomorphisms on a compact manifold, generated by solutions of
stochastic differential equations driven by a finite-dimensional
Brownian motion.

We show that {\it the presence of non-zero exponents combined with
certain nondegeneracy conditions} (amounting roughly speaking to
the assumption that the noise can move the orbit in any direction)
{\it implies almost surely chaotic behavior} in the following sense:
\begin{itemize}
\item{Exponential, in time, decay of  correlations between the
trajectories with different initial data.}
\item{Equidistribution of images of submanifolds.}
\item{ Central Limit Theorem, with respect to the measure on a
``rich enough'' subset, which holds for almost every fixed
realization  of the underlying Brownian motion.}
\end{itemize}

In order to illustrate the last point, let us consider a periodic
flow on $\R^n$, and let $\nu$ be a Lebegue probability measure
concentrated on an open subset. As a motivating example one may
think of an oil spot on the surface of the ocean.  The ultimate
goal could be to remove the oil  or at least to prevent it from
touching certain places on the surface of ocean. Thus, we wish to
predict the properties and the probability laws governing the
dynamics of the spot in time.
 Let $\nu_t$ be the measure on $\R^n$ induced from $\nu$ by
time $t$ map of the flow.  We shall show that almost surely $\nu_t$
is asymptotically equivalent to a Gaussian measure with variance of
order $t$. In other words, for a sufficiently large positive $R$
for large time 99 percent of the oil spot is contained in the ball
of radius  $R\sqrt t$.

Even though we consider the random flows  generated by  SDEs, very
little in our approach relies on the precise form of the noise,
and in a future work we shall show how to generalize our results
to other random flows.  Thus our work could be considered as a
first step in extending deterministic dynamical system picture to
a more general setup.

As a next step one may attempt to obtain the same results for
the so-called Isotropic Brownian Flows  introduced by Ito
(1956) and Yaglom (1957). This is a class of flows for which the
image of any simple point is a Brownian motion and the dependence
(the covariance tensor) between two different points is a function
of distance between these points. Related problems for
this case have been studied by
Harris \cite{H}, Baxendale \cite{BH}, Le Jan \cite{L1}, Cranston,
Scheutzow, Steinsaltz \cite{CSS1} and \cite{CSS2}, Cranston,
Scheutzow \cite{CS}, Lisei, Scheutzow \cite{LS},
Scheutzow, Steinsaltz \cite{SS} and others.

The precise statements of our results are given in the next
section. The proofs are carried out in Sections
\ref{Sc2PtCLT}--\ref{ScPathCLT}. Section \ref{ScDiss} deals
with dissipative flows. Application to passive transport
problems \cite{CC, CSS1} will be given elsewhere.

\section{Central Limit Theorems and an application to periodic flows.}

\subsection{Measure-Preserving Nondegenerate
Stochastic Flows of Diffeomorphisms}

Let $M$ be a $C^\infty$ smooth, connected, compact Riemannian
manifold with a smooth Riemannian metric $d$ and an associated
smooth measure $\mu$. Consider on $M$ a stochastic flow of
diffeomorphisms
\begin{equation}
\label{SF}
dx_t=\sum_{k=1}^d X_k(x_t) \circ d\theta_k(t)+ X_0(x_t) dt
\end{equation}
where $X_0, X_1, \dots , X_d$ are $C^\infty$-vector fields on $M$
and $\vec\theta(t)=(\th_1(t),\dots,\th_d(t))$ is a standard
$\R^d$-valued Brownian motion.  Since the differentials are in the
sense of Stratonovich, the associated generator for the process is
given by
\beq \label{gener}
L=\frac 12 \sum_{k=1}^d X_k^2 + X_0
\eneq

Let's impose additional assumptions on the vector fields
$X_0, X_1, \dots , X_d$. All together we impose five assumptions (A)
through (E). All these assumptions except (A) (measure preservation)
are {\it nondegeneracy
assumptions} and are satisfied for a generic set of vector fields
$X_0, X_1, \dots , X_d$. Now we formulate them precisely.

(A) ({\it measure preservation}) The stochastic flow $\{x_t: t\geq
0\}$, defined by (\ref{SF}), almost surely w.r.t. to the Wiener
measure $\cW$ of the Brownian motion $\vec\th$ preserves the
measure $\mu$ on $M$;

(B) ({\it hypoellipticity for $x_t$}) For all $x\in M$ we have
\beq
Lie(X_1,\dots,X_d)(x)=T_xM,
\eneq
i.e. the linear space spanned by all possible Lie brackets made out
of $X_1,\dots,X_d$ coincides with the whole tangent space $T_xM$
at $x$;

Denote by
\beq
\Delta=\{(x^1,x^2)\in M\times M:\ x^1=x^2\}
\eneq
the diagonal in the product $M \times M$.

(C$_2$) ({\it hypoellipticity for the two--point motion}) The
generator of the two--point motion $\{(x^1_t,x^2_t):\ t>0\}$ is
nondegenerate away from the diagonal $\Delta(M)$, meaning that the
Lie brackets made out of $(X_1(x^1),X_1(x^2)), \dots,
(X_d(x^1),X_d(x^2))$ generate $T_{x^1}M\times T_{x^2}M$.

To formulate the next assumption we need additional notations.
For $(t,x)\in [0,\infty) \times M$ let $Dx_t:T_{x_0}M \to T_{x_t}M$
be the linearization of $x_t$ at $t$. We need an analog of
hypoellipticity condition (B) for the process $\{(x_t,Dx_t):\ t>0\}$.
Denote by $TX_k$ the derivative of the vector field $X_k$ thought as
the map on $TM$ and by $SM=\{v\in TM:\ |v|=1\}$ the unit tangent
bundle on $M$. If we denote by $\~X_k(v)$ the projection of $TX_k(v)$
onto $T_vSM$, then the stochastic flow (\ref{SF}) on $M$ induces by
a stochastic flow on the unit tangent bundle $SM$ is defined by
the following equation:
\beq
d\~x_t(v)=\sum_{k=1}^d \~ X_k(\~x_t(v))\circ d\th_k(t)+
\~X_0(\~x_t(v))dt\ \ \textup{with}\ \~x_0(v)=v,
\eneq
where $v\in SM$. With these notations we have condition

(D) ({\it hypoellipticity for $(x_t,Dx_t)$}) For all $v\in SM$
we have
\[
Lie(\~X_1,\dots,\~X_d)(v)=T_vSM~.
\]

The next condition is that the stochastic flow $(x_t,Dx_t)$ has
a nonzero Lyapunov exponent. This condition can be expressed
in terms of a certain integral being nonzero. To define such
an integral precisely we need to introduce some additional objects.
Under condition (D) of hypoellipticity of $(x_t,Dx_t)$ there a unique
stationary probability measure $\~\mu$ on $SM$ for the one--point
motion $\{\~x_t:\ t>0\}$ on $SM$. Moreover, $\~\mu$ admits positive
density w.r.t. Riemannian measure on $SM$. Denote by $\pi:TM \to M$
the natural projection of $TM$ onto $M$. Consider the Levy-Civita
connection $\nabla$ for
the Riemannian structure on $M$ to define the horizontal space
$T_v TM$ for $v\in TM$ as the one identified with $T_{\pi v}M$.
Then $TX_k$ is the vector field on $TM$ for which $X_k(x)$ and
$\nabla X_k(x)(v)$ are the horizontal and vertical components of
$TX_k(v),\ v\in T_xM$ respectively. Define now the set of
$d+1$ functions
\beq
\beal
g_k(v)=<\nabla X_k(\pi v)(v),v>, \quad v\in SM,&\ \ k=0,1,\dots,d.
\\ R(v)=g_0(v)+\sum_{k=1}^d (L_{\~X_k} g_k)(v). & \
\enal
\eneq

With the above notations we have a formula of Carverhill \cite{C1}
for the largest Lyapunov exponent: \beq \lb_1=\lim_{t\to \infty}
\frac 1t \log |Dx_t(x)(v)|= \int_{SM} R(v) d\~\mu(v). \eneq (E)
({\it positive Lyapunov exponent}) Our last assumption is \beq
\lb_1>0 . \eneq For measure-preserving stochastic flows with
conditions (D) Lyapunov exponents $\lb_1,\dots, \lb_{\dim M}$ do
exist by {\it multiplicative ergodic theorem for stochastic flows}
of diffeomorphisms (see \cite{C2}, thm. 2.1). Moreover, the sum of
Lyapunov exponents $\sum_{j=1}^{\dim M}\lb_j$ should be zero (see
e.g. \cite{BS}). Therefore, to have a positive Lyapunov exponent
it is sufficient to have a nonzero Lyapunov exponent. Lyapunov
exponents of the stochastic flows of diffeomorphisms have been
investigated by various authors, see \cite{Ba, C1, ES} and
references there.

\subsection{CLTs for Measure-Preserving Nondegenerate
Stochastic Flows of Diffeomorphisms}

From now on we shall consider only the stochastic flows, defined
by (\ref{SF}), and satisfying assumptions (A) through (E). The
first CLT for such flows is CLT for additive functionals of the
two-point motion.

Denote by $\{A_t^{(2)}:\ t>0\}$ an additive functional of
the two--point motion $\{(x^1_t,x^2_t):\ t>0\}$. Suppose $A_t^{(2)}$
is governed by the Stratonovich stochastic differential equation
\begin{equation}
\label{Add2}
dA_t^{(2)}(x^1_0,x^2_0)=
\sum_{k=1}^d \alpha_k(x^1_t, x^2_t)\circ d\theta_k(t)
+ a(x^1_t,x^2_t) dt,
\end{equation}
where $\{\alpha_k\}_{k=1}^d$ and $a$ are $C^\infty$-smooth
functions, such that
\beq \label{2pt}
\beal
\int\!\!\!\!\int_{M\times M} \Big[a(x^1,x^2) +& \\
\frac{1}{2}\sum_k
\left(L_{(X_k,X_k)}\alpha_k\Big)(x^1,x^2)\right] & \
d\mu(x^1)\ d\mu(x^2)=0
\enal
\eneq
and $\left(L_{(X_k,X_k)}\alpha_k\right)(x^1,x^2)$
denotes the derivative of $\alpha_k$ along the vector field
$(X_k,X_k)$ on $M \times M$ at the point $(x^1,x^2)$.
This equality can be attained by subtracting a constant from
$a(x^1,x^2)$.

\bthm \label{CLT2pt} Let $\{A_t^{(2)}:\ t>0\}$ be the  additive
functional of the  two--point motion, defined by (\ref{Add2}), and
let the  two--point motion $\{(x^1_t,x^2_t):\ t>0\}$ for
$x^1_0\neq x^2_0$ be defined by the stochastic flow of diffeomorphisms
(\ref{SF}), which in turn satisfies  conditions (A) through (E).
Then as $t\to\infty$ we have that $\frac{A_t^{(2)}}{\sqrt t}$
converges weakly to a normal random variable.
\ethm

Fix a positive integer $n> 2$. The second CLT for stochastic
flows, defined by (\ref{SF}) and satisfying assumptions (A)
through (E), is the CLT for additive functionals of the  $n$-point
motion. In other words, CLT for the  two--point motion (Theorem
\ref{CLT2pt}) has a natural generalization to a CLT for the
$n$-point motion. Denote by \beq \nonumber
\Delta^{(n)}(M)=\{(x^1,\dots,x^n)\in M \times \dots \times M: \
\exists j\neq i \ \ \textup{such that}\ \  x^j=x^i\} \eneq the
generalized diagonal in the product $M \times \dots \times M$ ($n$
times). Replace (C$_2$) by the following condition.

(C$_n$) The generator of the  $n$--point motion
$\{(x^1_t,\dots,x^n_t):\ t>0\}$ is non-degenerate away from the
generalized diagonal $\Delta^{(n)}(M)$, meaning that Lie brackets
made out of $(X_1(x^1),\ \dots,\ X_1(x^n)),\ \dots,\ (X_d(x^1),\
\dots,\ X_d(x^n))$ generate $T_{x^1}M\times \dots \times
T_{x^n}M$.

Similarly to the case of the  two--point motion, denote by
$\{A^{(n)}_t:\ t>0\}$ an additive functional of the  $n$-point
motion $\{(x^1_t,\dots,x^n_t):\ t>0\}$. Suppose $A^{(n)}_t$ is
governed by the Stratonovich stochastic differential equation
\beq
\label{Addn} dA^{(n)}_t(x^1_0, \dots x^n_0)= \sum_{k=1}^d
\alpha_k(x^1_t, \dots, x^n_t)\circ d\theta_k(t)+ a(x^1_t,\dots,
x^n_t) dt,
\eneq
where $\{\alpha_k\}_{k=1}^d$ and $a$ are
$C^\infty$-smooth functions, which satisfy
\beq
\beal \label{npt}
\int \cdots
\int_{M\times \dots \times M} \Big[a(x^1,\dots,x^n)+ & \\
\frac{1}{2}\sum_{k=1}^d \left(
L_{(X_k,\dots,X_k)}\alpha_k\Big)(x^1,\dots,x^n) \right] & \
d\mu(x^1) \dots d\mu(x^n)=0~,
\enal
\eneq
and $\left(L_{(X_k,\dots,X_k)}\alpha_k\right)(x^1,\dots,x^n)$
denotes the derivative of $\alpha_k$ along the vector field
$(X_k,\dots,X_k)$ on $M \times \dots \times M$ ($n$ times) at the
point $(x^1,\dots,x^n)$. As in the two--point case, this equality
can be attained by subtracting a constant from $a(x^1,\dots,x^n)$.
For the $n$-point motion we have the following

\begin{theorem} \label{CLTnpt} Let $\{A_t^{(n)}:\ t>0\}$ be
the additive functional of the $n$-point motion, defined by
(\ref{Addn}), and let the $n$-point motion
$\{(x^1_t,\dots,x^n_t):\ t>0\}$ for pairwise distinct
$x^i_0\neq x^j_0$ be defined by the stochastic flow of
diffeomorphisms (\ref{SF}), which in turn satisfies conditions
(A),(B), (C$_n$), (D), and (E). Then as $t\to\infty$ we have
that $\frac{A^{(n)}_t}{\sqrt t}$ converges weakly to a normal
random variable.
\end{theorem}
We shall prove the CLT for the  two--point motion (Theorem
\ref{CLT2pt}), and then show how to extend this proof to the
$n$-point case.

The third CLT for stochastic flows, defined by (\ref{SF}) and
satisfying assumptions (A) through (E), is the CLT for probability
measures supported on sets of positive Hausdorff dimension in $M$.

Consider stochastic flow (\ref{SF}) and an additive functional
$\{A_t:\ t>0\}$ of the  one--point motion satisfying \beq
\label{Add1} dA_t(x)=\ \sum_{k=1}^d \alpha_k(x_t) \circ
d\theta_k(t)+\ a(x_t) dt \eneq with $C^\infty$-smooth
coefficients. Define \beq \ha(x)=\
\sum_{k=1}^d(L_{X_k}\alpha_k)(x)+\ a(x), \eneq where
$(L_{X_k}\alpha_k)(x)$ is the derivative of $\alpha_k$ along the
vector field $X_k$ at point $x$. Impose additional assumptions on
the coefficients of (\ref{Add1}).

(F) ({\it no drift or preservation of the center of mass})
\begin{equation} \label{CenterOfMass}
\beal \int_M \ha(x) d\mu(x)=0~, \quad \int_M \alpha_k(x) d\mu(x)=0
\quad \textup{for} \quad k=1,\dots,d. \enal
\end{equation}
This condition can be attained by subtracting appropriate constants
from functions $\alpha_1,\dots,\alpha_d,$ and $a$. For $A_t$ as above,
when $t\to\infty$, we have that $\frac{A_t}{\sqrt t}$ converges to
a normal random variable with zero mean and some variance $D(A)$,
given below. Our next result below shows how little randomness in
initial condition is needed for the CLT to hold.

Let $\nu$ be a probability measure on $M$, such that for some
positive $p$ it has a finite $p$-energy
\beq \label{energy}
I_p(\nu)=\int\!\!\!\!\int\frac{d\nu(x) d\nu(y)}{d^p(x,y)}<\infty.
\eneq
In particular, this means that the Hausdorff dimension of
the support of $\nu$ on $M$ is positive (see \cite{Ma}, sect. 8).
Let $\cM^\theta_t$ be the measure on $\R$ defined on Borel sets
$\Omega \subset \R$ by
\beq \label{indmeas}
\cM^\theta_t(\Omega)=\nu\left\{x\in M:\
\frac{A^\theta_t(x)}{\sqrt t}\in\Omega\right\}.
\eneq
\begin{theorem} \label{measure}
Let $\{x_t:\ t>0\}$ be a stochastic flow of diffeomorphisms
(\ref{SF}), and let  conditions (A) through (F) be satisfied. Then
as $t\to+\infty$  almost surely $\cM_t^\theta$ converges weakly to
the Gaussian measure with zero mean and some variance $D(A)$.
\end{theorem}

\subsection{Application to periodic flows.}

Consider the stochastic flow (\ref{SF}) on $\R^N$, with the periodic
vector fields $X_k$. Application of Theorems \ref{CLT2pt}--\ref{measure}
to the corresponding flow  on the $N$-dimensional torus leads to
the clear statements on the behavior of the flow on $\R^N$. We formulate
those as Theorems 1$'$ -- 3$'$ below.

The usual CLT describes the distribution of the displacement of a single
particle with respect to the measure of the underlying Brownian motion
$\theta(t)$. The CLT formulated below (Theorem 3$'$), on the contrary,
holds for almost every realization of the Brownian motion and is with
respect to the randomness in the initial condition.

{\bf Theorem 1$'$.}\ {\it Let $\{X_k\}_{k=0}^d$ be $C^{\infty}$
periodic vector fields in $\R^N$ with a common period, and let
conditions (A) through (E) be satisfied. Let $x^1_t$ and $x^2_t$
be the solutions of (\ref{SF}) with different initial data. Then
for some value of the drift $v$  the vector $\frac{1}{\sqrt{t}}
(x^1_t-v t, x^2_t-v t)$ converges as $t \to \infty$ to a Gaussian
random vector with zero mean.}

{\bf Theorem 2$'$.}\ {\it Let $\{X_k\}_{k=0}^d$ be $C^{\infty}$
periodic vector fields in $\R^N$ with a common period, and let the
conditions $(A),(B),(C_n),(D),$ and $(E)$ be satisfied. Let
$x^1_t,\dots, x^n_t$ be the solutions of (\ref{SF}) with pairwise
different initial data. Then for some value of the drift $v$ the
vector $\frac{1}{\sqrt{t}} (x^1_t-vt,..., x^n_t-vt)$ converges as
$t \to \infty$ to a Gaussian random vector with zero mean.}

{\bf Theorem 3$'$.}\ {\it Let $\{X_k\}_{k=0}^d$ be $C^{\infty}$
periodic vector fields in $\R^N$ with a common period, let $\nu$
be a probability measure with finite $p$-energy for some $p>0$ and
with compact support, and let conditions (A) through (F) be
satisfied. For the condition (F) we take $\alpha_k = X_k$ and $a =
X_0$. Let $x_t$ be the solution of (\ref{SF}) with the initial
measure $\nu$. Then for almost every realization of the Brownian
motion the distribution of $\frac{x_t}{ \sqrt{t}}$ induced by
$\nu$ converges weakly as $t \to \infty$ to a Gaussian random
variable on  $\R^N$ with zero mean and some variance $D$.}

{\it Proof of Theorems 1$'$ - 3$'$:}\ The functions $\alpha_k$
and $a$ in the formulas (\ref{Add2}), (\ref{Addn}), and (\ref{Add1})
could be considered to be vector valued, thus defining the vector
valued additive functionals. Any linear combination of the
components of a vector valued additive functional is a scalar
additive functional, for which Theorems
\ref{CLT2pt}--\ref{measure} hold. If any linear combination of
the components of a vector is a Gaussian random variable then
the vector itself is Gaussian. Therefore Theorems
\ref{CLT2pt}--\ref{measure} hold for vector valued additive
functionals as well.

It remains to rewrite equation (\ref{SF}) in the integral form and
apply Theorems \ref{CLT2pt}, \ref{CLTnpt}, or \ref{measure} to
the vector valued additive functional of the flow on the torus. Q.E.D.

\begin{remark}
Let in Theorem 2$'$
$$\frac{1}{\sqrt{t}} (x^1_t-vt,\dots, x^n_t-vt) \to
\left(\Theta^1,\dots, \Theta^n\right)$$
where $\Theta^j\in\reals^N$ are Brownian motions. If $x_t$ moves
according to (\ref{SF}) and satisfies condition (F), then $\Theta^j$
and $\Theta^i$ are independent for $j\neq i.$ In fact it will follow
from the proof of Theorem 2 that components of $\Theta$ satisfy
\beq
\beal
\EXP\left(\Theta^j_l \Theta^i_r\right) & =
\lim_{t\to\infty} \frac{1}{t}\int\!\!\!\!\int\left[
\left(\int_0^t (\hX_0)_l (x^j_s) d (x^i_s)_r\right)
d\mu(x^j_0) d\mu (x^i_0)+ \right. \\
& \left.\left(\int_0^t (\hX_0)_r (x^i_s) d (x^j_s)_l\right)
d\mu(x^j_0) d\mu (x^i_0)+ \right. \\
& \left.
\sum_{k=1}^d \left((X_k)_l (x^j) (X_k)_r (x^i)\right) d s\right]
d\mu(x^j) d\mu(x^i),
\enal
\eneq
where $\hX$ is the lift to the universal cover $\R^N$ of
the vector field $X$, $\Theta^j_l$ is the $l$-th component of
$\Theta^j$, and $x^j_s$ is the $s$-th component of $x^j$
(cf. Lm. \ref{L5}). However performing the integration over
$\mu\times\mu$ first we get that all terms equal 0.
\end{remark}

The proofs of the CLT's occupy Sections \ref{Sc2PtCLT}--\ref{ScPathCLT}.
In the next section
we prove CLT for additive functionals of the two--point motion
(Theorem \ref{CLT2pt}). By results of Baxendale and Stroock we
know that the  two--point motion is ergodic \cite{BS}. In section
\ref{mixing} we investigate the rate of mixing (the decay of
correlations) of the two--point process, and prove that it is
exponential in time. In the next section we show how, using the
proof for two--point functionals, one can derive CLT for the
$n$-point functionals (Theorem \ref{CLTnpt}). In section
\ref{curve} we prove that a smooth curve $\gamma$ on $M$ becomes
uniformly distributed by flow (\ref{SF}) in the limit as time tends
to infinity. Section \ref{ScPathCLT} of this paper is devoted to the
proof of CLT for measures (Theorem \ref{measure}).
In Section \ref{ScDiss} we prove above CLT for the dissipative  case.

\section{Proof of CLT for two--point functionals}
\label{Sc2PtCLT}
\subsection{Outline of the proof}

We consider the  two--point motion $z_t=(x^1_t,x^2_t)$, defined by
(\ref{SF}). The generator of the two--point process degenerates on
the diagonal $\Delta \subset M \times M$, and this creates a
problem when we want to establish the mixing properties of the
two--point process, needed for the CLT. However, since we consider
stochastic flows (\ref{SF})  with positive Lyapunov exponents, the
intervals of time when the points $x^1_s$ and $x^2_s$ are nearby
should be small. This
observation will help us to detour the problem of the diagonal.

The outline of the proof of the CLT is the following. Choose a
small $r > 0$ , and two initial points $x^1\neq x^2$. Denote
$z=(x^1,x^2)$ and $z_t=(x^1_t, x^2_t).$ Let
\beq \label{r-nbhd}
\beal
\Delta_r=\{(x,y)\in M\times M: d(x,y)=r\},\\
G_r=\{(x,y)\in M\times M: d(x,y)>  r\}.
\enal
\eneq

We introduce two sequences of stopping times $\sigma_n$ and $\tau_n$,
such that $z_{\sigma_n} \in \Delta_{r/2}$ and $z_{\tau_n} \in G_r$.
These stopping times are defined inductively, with $\sigma_n$ being
the first instance after $\tau_{n-1}$ when the process $z_t$ visits
$\Delta_{r/2}$, and $\tau_n$ related to the first instance after
$\sigma_n$ when the process $z_t$ visits $\Delta_r$. We shall see
that $z_{\tau_n}$ is an exponentially mixing Markov chain, and that
the stopping times $\tau_n$ satisfy the law of large numbers,
$\tau_n /n \rightarrow {\rm const}$. The contributions to
the additive functional $A_t$ from the intervals
$(\tau_{n-1}, \tau_n)$ will be seen to satisfy the CLT.

Up to some technical details this is the main idea of the proof.
In the proof we shall use results of Baxendale--Stroock \cite{BS}
on the properties of the  stopping times of the two--point motion
of stochastic flows (\ref{SF}) with nonzero Lyapunov exponents,
and standard facts about hypoelliptic operators.

\subsection{Preparatory Lemmas and the proof}
\label{pl}

In this section we prove several Lemmas leading to the CLT for the
two--point motion. First, we study the two--point motion close to
the diagonal $\Delta \subset M \times M$. We shall use the
following result of Baxendale--Stroock \cite{BS} showing, in
particular, that a positive Lyapunov exponent for (\ref{SF})
implies that the transition time from $r/2$ neighborhood of the
diagonal to $r$ neighborhood of the diagonal has exponential
moments.

\begin{lemma} \label{Escape} (\cite{BS} see Theorem 3.19)
There are positive constants $\alpha_0$ and $r_0$, which depend
on the vector fields $X_0,\dots , X_d$, such that for any
$0<\alpha<\alpha_0$ and $0<r<r_0$, there are $p(\alpha),\
K(\alpha)$ with the property $p(\alpha)\to 0$ as $\alpha\to 0$,
such that  for any pair of distinct points $x^1$ and $x^2$ which
are at most $r$ apart we have
\beq \label{ine1}
K^{-1}\
d^{-p}(x^1,x^2)\leq E_z( e^{\alpha \tau})\leq K\ d^{-p}(x^1,x^2),
\eneq
where $\tau$ is the stopping time of $x^1_s$ and $x^2_s$
getting distance $r$ apart, i.e. $\tau=\inf\{s>0:\
d(x^1_s,x^2_s)=r\}$ and $z=(x_0^1,x_0^2)$.
\end{lemma}

Fix some $r>0$, which satisfies the assumptions of
 Lemma \ref{Escape}. Let $\delta$ be small enough, so that
\begin{equation} \label{formula1}
\PROB_{z_0} \{z_t \in \Delta_{r/2}~~~{\rm for}~{\rm some}~~0\leq t
\leq 2\delta \} \leq \frac{1}{2K}~,
\end{equation}
for all $z_0 \in \Delta_r$, where the constant $K$ is the same as
in (\ref{ine1}). Let $\Z_{\delta}^+ = \{ \delta, 2\delta,... \}$.
Next we introduce the stopping times $\sigma_n, \tau_n,~n\geq 1$.
Let $\sigma_1 = \inf \{t: z_t \in \Delta_{r/2} \}$. Assuming that
$\sigma_n$ has been defined, we define $\tau_n$ as follows
\beq \label{stop}
\beal
\tau_n = \inf \{t> \sigma_n: t \in \Z_{\delta}^+, {\rm there}~{\rm
exists}~ s_1 \in [\sigma_n , t-\delta],~{\rm such}~{\rm that}~ \\
z_{s_1} \in \Delta_r, ~
z_s \in G_{r/2}~{\rm  for}~{\rm all }~ s \in [s_1, t] \}
\enal
\eneq
This definition means that if a trajectory starts at a point in
$\Delta_{r/2}$, then we wait till it hits the set $\Delta_r$, and
stop in the second of the discrete moments of time, following
the visit to $\Delta_r$, provided the part of the trajectory after
the time it hits $\Delta_r$ is entirely contained in $G_{r/2}$.
Finally, if $\tau_n$ is defined, we define
$\sigma_{n+1} = \inf \{ t > \tau_n: z_t \in \Delta_{r/2} \} $.

We have constructed the Markov chain  $z_{\tau_n}$, with
the state space $G_{r/2}$. With the  discretization of time
we were able to stop the process $z_t$ in the open set, and will
now show that, due to hypoellipticity, the Markov chain $z_{\tau_n}$
is exponentially mixing.

\begin{lemma}
\label{Induced} The Markov chain  $z_{\tau_n}$ is exponentially
mixing.
\end{lemma}
{\it Proof:} Let $U$ be some open set, which, together with its
closure, is contained in $G_{r/2}$. Let $P(x,y)$ be the density of
the absolutely continuous component  (with respect to the measure
$\mu^{(2)}$, invariant for the two--point process) of the
transition function for the chain $z_{\tau_n}$\footnote{It is
shown in \cite{K} that under the condition (A) the invariant
measure for the two--point process is the product measure, i.e.
$\mu^{(2)}=\mu \times \mu$}. It is sufficient to show (see
\cite{Do}), that $P(x,y) \geq c$ for all $x \in G_{r/2}~$, $y \in
U$, and some $c >0$. Consider the process in $\overline{G}_{r/2}$,
which coincides with $z_t$ in $G_{r/2}$, and is stopped, when it
reaches the boundary. Let $p_t(x,y),~~ x,y \in \overline{G}_{r/2}
$ be the density of the stochastic transition function for this
process.

 Let $K$ be a compact in $G_{r/2}$, which contains both
$U$ and $\Delta_r$. Due to hypoellipticity for  each point $x \in
G_{r/2}$ there is a neighborhood $V(x)$, and some $T(x) > 0$, such
that $p_t(y_1, y_2) > c(x, t)>0$ when $y_1, y_2 \in V(x)$ and $t
\in (0, T(x)]$  (see \cite{Le}). The lower bound $c(x, t)$ is
strictly positive on time intervals separated from zero $c(x,t) >
c(x) > 0$ for $t \in [\varepsilon, T(x)]$ for any positive
$\varepsilon$. By covering the compact $K$ with a finite number of
such neighborhoods we see that the density is estimated from below
\[
p_t(x,y) \geq c > 0~, ~~~~~~~~{\rm for}~~{\rm all}~ x \in
\Delta_r, ~y \in U,~ t \in [\delta, 2\delta]~.
\]
From the definition of the stopping time $\tau_1$ if follows that
\beq \label{estBel}
P(x,y) \geq \inf_{x \in \Delta_r, ~ t \in
[\delta, 2\delta]} p_t(x,y) >c.
\eneq
This implies the Lemma. Q.E.D.

The next Lemma shows that the distributions of the transitions
times have uniformly exponentially decreasing tails.

\begin{lemma} \label{RetExpMom}  There exist positive $D$ and
$\gamma$, such that for all $z \in G_{r/2}$ we have
\beq
\PROB_z \{\tau_1>t \}\leq D\ e^{-\gamma t}.
\eneq
\end{lemma}
{\it Proof:} Let us introduce two more sequences of stopping times,
$a_n$ and $b_n$, $n \geq 1$. These are the times of the successive
visits by the process of the sets $\Delta_{r/2}$ and $\Delta_r$
respectively:
\[
a_1 = \inf\{t: z_t\in \Delta_{r/2}\},~~ b_n = \inf\{t > a_n:
z_t\in \Delta_r \},
\]
\[
  a_{n+1} = \inf\{t > b_n: z_t\in
\Delta_{r/2} \}~.
\]
Let $E_n$ be the event that the process visits $\Delta_{r/2}$, and
then  makes $n$ excursions before time $\tau_1$ from $\Delta_r$ to
$\Delta_{r/2}$ and back: $ E_n = \{ b_{n+1} < \tau_1 < b_{n+2}
\}$. Then,
\[ \EXP_z
e^{\alpha \tau_1} = \sum_{n=0}^{\infty} \EXP_z (e^{\alpha \tau_1}
\chi_{E_n} )~. \]
Notice that $\tau_1 \leq b_{n+1} + 2\delta$ on
the set $E_n$, and the set $E_n$ is contained in the set
$\{a_{i+1} - b_i \leq 2 \delta~~~{\rm for}~{\rm all}~~1 \leq i
\leq n \}$. Thus,
\[
\EXP_z e^{\alpha \tau_1} \leq \sum_{n=0}^{\infty}\EXP_z (
e^{\alpha(b_{n+1} + 2\delta)}\chi_{\{a_{i+1}-b_i \leq
2\delta~~{\rm for}~~{\rm all }~~1 \leq i
\leq  n\}})~.
\]
By the Markov property, for $n \geq 1$ we have
\[
\EXP_z \left(
e^{\alpha(b_{n+1} + 2\delta)}
\chi_{\{a_{i+1}-b_i \leq 2\delta~~{\rm for}~~{\rm all }~~1
\leq i \leq  n\}} \right) =
\]
\[
\EXP_z \left( e^{\alpha(b_{n} +
2\delta)}
\chi_{\{a_{i+1}-b_i \leq 2\delta~~{\rm for}~~{\rm all }~~
1 \leq i \leq  n-1\}}\ \right)
\EXP_{z_{b_{n}}}( e^{\alpha b_1} \chi_{\{a_1 \leq 2\delta \} })~.
\]
For the next estimate we use the result of Lemma \ref{Escape},
together with (\ref{formula1}),
\[
\sup_{z \in \Delta_r} \EXP_{z}( e^{\alpha b_1} \chi_{\{a_1 \leq
2\delta \} }) \leq e^{2\delta \alpha} \frac{1}{2K} \sup_{z \in
\Delta_{r/2}} \EXP_z e^{\alpha \tau} \leq 1/2 e^{2 \delta \alpha}
\left(\frac{2}{r}\right)^{p(\alpha)}~.
\]
This can be made smaller than $3/4$ by selecting $\alpha$
sufficiently small. Therefore,
\beq
\beal
\EXP_z \left( e^{\alpha(b_{n+1} + 2\delta)} \chi_{\{a_{i+1}-b_i \leq
2\delta~~{\rm for}~~{\rm all }~~1 \leq i \leq  n\}} \right)
\leq \\
\frac 34
\ \EXP_z \left( e^{\alpha(b_{n} + 2\delta)} \chi_{\{a_{i+1}-b_i \leq
2\delta~~{\rm for}~~{\rm all }~~1 \leq i \leq  n-1\}}~ \right).
\enal
\eneq
Continuing by induction we see that
\[
\EXP_z \left(e^{\alpha(b_{n+1} +
2\delta)}
\chi_{\{a_{i+1}-b_i \leq 2\delta~~{\rm for}~~{\rm all }~~1 \leq i
< n\}}\right)
\leq
\left(\frac{3}{4}\right)^n  \EXP_z e^{\alpha(b_{1} + 2\delta)}~.
\]
Thus,
\[
\sup_{z \in G_{r/2}}  \EXP_z e^{\alpha \tau_1} \leq  \sup_{z \in
G_{r/2}} \EXP_z e^{\alpha(b_{1} + 2\delta)} \sum_{n=0}^{\infty}
(3/4)^n = D~.
\]
This implies the statement of the Lemma. Q.E.D.

Consider now the functional $A_t^{(2)}(x^1,x^2)$ of the two point
motion, defined by the stochastic equation (\ref{Add2}). Using
Lemma \ref{RetExpMom} we now show the existence of exponential
moments for $A^{(2)}_{\tau_1}(z)$
\begin{lemma} \label{RetExpMom2}  There exist positive $D_1$ and
$\gamma_1$, such that for all $z \in G_{r/2}$ we have
\beq
\PROB_z \{|A^{(2)}_{\tau_1}|>t \}\leq D_1\ e^{-\gamma_1 t}.
\eneq
\end{lemma}
{\it Proof:}  Rewriting the Stratonovich equation
(\ref{Add2}) in terms of Ito equation we get \beq \beal
\label{Ito} A^{(2)}_{\tau_1}(z)=\sum_{k=1}^d \int_{0}^{\tau_1}
\alpha_k(z_s) d\theta_k(s) +
\\ \int_{0}^{\tau_1} \left[a(z_s)+
\frac{1}{2}\sum_{k=1}^d (L_{X_k,X_k}\alpha_k)(z_s)\right] ds, \enal
\eneq where $(L_{X_k,X_k}\alpha_k)(z_s)$ is the derivative of
$\alpha_k$ along the vector field $(X_k,X_k)$ on $M \times M$ at
the point $z_s$. Let
\beq
\beal \label{Bs}
B_0(t) =  \int_{0}^{t} \left[a(z_s)+
\frac{1}{2}\sum_{k=1}^d (L_{X_k,X_k}\alpha_k)(z_s)\right] ds~,~~{\rm
and}
\\ B_k(t) = \int_{0}^{t} \alpha_k(z_s) d\theta_k(s),
~~k=1,\dots,d~.
\enal
\eneq
In order to prove the Lemma it is sufficient to show that
\[
\EXP_z e^{\delta | B_k(\tau_1)| } < \infty~,~~~k = 0,...,d
\]
for a sufficiently small $\delta$. By Lemma \ref{RetExpMom} we
have
\beq \label{eq3}
\beal \EXP_z e^{\delta | B_k(\tau_1)| } \leq
& \sum_{n = 0}^{\infty}\PROB_z\{n < \tau_1 \leq n+1\}\  \EXP_z
(e^{\delta | B_k(\tau_1)| } \chi_{\{\tau_1 \leq n+1\}}) \leq \\ &
\sum_{n = 0}^{\infty} D e^{-\gamma n } \ \EXP_z e^{\delta
\sup_{[0,n+1]} |B_k(t)|}
\enal
\eneq
Notice that for $k=1,\dots ,d$
the integral $B_k(t)$ is obtained from a Brownian motion by a
random time change. This time change has a bounded derivative,
since the functions $\alpha_k$ are bounded. If $w(t)$ is a
Brownian motion, then for any $\alpha $ there is $q(\alpha)$ with
the property $q(\alpha) \rightarrow 0 $ when $ \alpha \rightarrow
0$, such that \beq \EXP (e^{\alpha \sup_{[0,c(n+1)]} |w(t)|})<C
e^{q(\alpha) n} . \eneq For $k =0$ we have $ |B_0(t)| < Ct$ since
the integrand in (\ref{Bs}) is bounded. Thus the right hand side
of (\ref{eq3}) is finite for a sufficiently small $\delta$ for all
$k$. This completes the proof of the Lemma. Q.E.D.

Let $\mu_0$ be the measure on $G_{r/2}$ which is invariant for
the Markov chain  $z_{\tau_n}$. Consider two Markov chains with
the state space $G_{r/2} \times \mathbb{R}$, defined by their
stochastic transition functions
\beq  \label{PP}
\beal
P_z(dy \times dt)
= & \PROB_z \{z_{\tau_1} \in dy; \tau_1 \in dt\}~, \\
Q_z(dy \times dt) = & \PROB_z \{z_{\tau_1} \in dy;
A^{(2)}_{\tau_1} \in dt\}~.
\enal
\eneq
Note that the transition functions are defined to depend only on
the first component of the original point.  The measures
\beq
\beal
\mu^{\tau} (dy \times dt) = & \ \int P_z \{dy \times dt\}\
d\mu_0(z)~,~~~~~{\rm and}
\\
\mu^{A} (dy \times dt) = & \ \int Q_z (dy \times dt)\ d \mu_0(z)
\enal \eneq are invariant for their respective Markov chains.
These invariant measures are unique, since $\mu_0$ is the unique
invariant measure for the process $z_{\tau_n}$. As easily follows
from (\ref{estBel}) both of these Markov chains satisfy the
Doeblin condition.

Notice that transition probabilities for the Markov chains
$(z_{\tau_n},\tau_n(z) - \tau_{n-1}(z))$ and
$( z_{\tau_n}, A^{(2)}_{\tau_n}(z) - A^{(2)}_{\tau_{n-1}}(z))$
are given by (\ref{PP}). For each of the Markov chains we apply
the Law of Large Numbers to the function
$f^{\tau (\textup{resp.}\ A)}: G_{r/2} \times \mathbb{R}
\rightarrow \mathbb{R}$ given by $f(z,t) = t$. Notice that
$f^\tau \in L_1(d\mu^{\tau})$ and $f^A \in L_1(d\mu^{A})$ by Lemmas
\ref{RetExpMom} and \ref{RetExpMom2}. Since both Markov chains are
ergodic, in the first case we obtain for each $z \in G_{r/2}$
\beq
\beal \label{taul}
\lim_{n \rightarrow \infty}
\frac{\tau_n(z)}{n} = \lim_{n \rightarrow \infty} \frac{\tau_1+
(\tau_2(z) - \tau_1(z))+...+( \tau_n(z)- \tau_{n-1}(z)) }{n}\\
= \int_{G_{r/2}} \EXP_z \tau_1 (z) d \mu_0 (z)
~~~~{\rm almost}~~{\rm surely}~,
\enal
\eneq
while for the second Markov chain we get for $z \in
G_{r/2}$
\beq  \label{Al}
\beal
\lim_{n \rightarrow \infty}
\frac{A^{(2)}_{\tau_n}(z)}{n} = \lim_{n \rightarrow \infty}
\frac{A^{(2)}_{\tau_1}+ (A^{(2)}_{\tau_2} - A^{(2)}_{\tau_1})+
\dots + ( A^{(2)}_{\tau_n}- A^{(2)}_{\tau_{n-1}}) }{n} = \\
\int_{G_{r/2}} \EXP_z A^{(2)}_{\tau_1} (z) d \mu_0 (z) ~,
~~~~{\rm almost}~~{\rm surely}~
\enal
\eneq
\begin{lemma} \label{lemmaA}
For any $z \in G_{r/2}$ the limit $\lim_{n \rightarrow \infty}
\frac{A^{(2)}_{\tau_n}(z)}{n}$ is equal to zero almost surely.
\end{lemma}
{\it Proof:}\ Since the Markov chain is ergodic, we only need to
show that the limit
$\lim_{n \to \infty} \frac{\EXP_z A^{(2)}_{\tau_n}}{n}$
is zero for any $z \in G_{r/2}$. As in (\ref{Ito}),  we have
\beq
\beal \label{Ito2} A^{(2)}_{\tau_n}(z)=\sum_{k=1}^d
\int_{0}^{\tau_n} \alpha_k(z_s) d\theta_k(s) +
\\
\int_{0}^{\tau_n} \left[a(z_s)+
\frac{1}{2}\sum_{k=1}^d (L_{X_k,X_k}\alpha_k)(z_s)\right] ds,
\enal
\eneq
Taking expectations of both sides we obtain
\[
\EXP_z A_{\tau_n}(z) = \EXP_z \int_{0}^{\tau_n} \varphi(z_s) ds~,
\]
where $\varphi(z) = a(z)+ \frac{1}{2}\sum_{k=1}^d
(L_{X_k,X_k}\alpha_k)(z)$ is a function, whose integral with
respect to the measure $\mu^{(2)}$, invariant for the two--point
process, is equal to zero by (\ref{2pt}). Since the two--point
motion is ergodic we have that
\[
\lim_{t \rightarrow \infty} \frac{\EXP_z \int_0^t \varphi(z_s)
ds}{t} = \int\!\!\!\!\int \varphi(z) d\mu^{(2)}(z) = 0~.
\]
Therefore,
\beq
\beal
\lim_{n \rightarrow \infty} \frac{\EXP_z \int_0^{\tau_n}
\varphi(z_s) ds}{n} = & \lim_{n \rightarrow \infty}
\frac{\tau_n}{n} \frac{\EXP_z \int_0^{\tau_n} \varphi(z_s) ds}{\tau_n}
\\
= \int_{G_{r/2}} \EXP_z \tau_1 (z) d \mu_0 (z) &
~\lim_{t \rightarrow \infty}
\frac{\EXP_z \int_0^t \varphi(z_s) ds}{t}=0~.
\enal
\eneq
This completes the proof of the Lemma. Q.E.D. \\

\brm We have shown that
\beq \label{zer} \int_{G_{r/2}}
\EXP_z A_{\tau_1} (z) d \mu_0 (z) = 0~.
\eneq
In section \ref{mixing} we shall use the identity
\[
\int_{G_{r/2}} \int_0^{\tau_1} \EXP_z  e^{i 2 \pi k \delta^{-1} s}
\varphi (z_s)  ds d \mu_0 (z) = 0
\]
for integer $k$ and a function $\varphi$, such that  $\int
\varphi(z) d \mu^{(2)}(z) = 0$. It can be proved similarly to the
way formula (\ref{zer}) was proved if we recall that $\tau_n \in
\{ \delta, 2\delta,\dots \} $. \erm

{\it Proof of Theorem \ref{CLT2pt}:}\ Let
\[
N(t) = \left[\frac{t}{\int \EXP_z \tau_1(z) d \mu_0(z)}\right]~.
\]
Then,
\beq \label{st}
\frac{A_t^{(2)}}{\sqrt{t}} =
\frac{A_{\tau_{N(t)}}^{(2)}}{\sqrt{t}} +
\frac{  A_t^{(2)} - A_{\tau_{N(t)}}^{(2)}}{\sqrt{t}}~.
\eneq
Without loss of generality we may consider the initial point $z$
distributed in $G_{r/2}$ according to some law. Lemma \ref{lemmaA}
implies that the function $f^A(z,t) = t$ on $ G_{r/2} \times \mathbb{R}$
has zero mean with respect to measure $\mu^A$. We can, therefore,
use the CLT for Markov chains to obtain
\beq \label{forml}
\lim_{t \rightarrow \infty} \frac{A^{(2)}_{\tau_{[s N(t)]}} }{\sqrt{t}}
= \Theta(s)~,
\eneq
where $\Theta(s)$ is a Wiener process with some
diffusion coefficient $\sigma$. In particular,
\beq
\frac{A^{(2)}_{\tau_{N(t)}}}{\sqrt{t}} \rightarrow N(0, \sigma)~,
\label{ab0} \eneq and for any $c > 0$ we have \beq \label{ab1}
\sup_{r \in [N(t) -c \sqrt{t}, N(t) + c \sqrt{t}] \bigcap
\mathbb{N}}
 \frac{|A^{(2)}_{\tau_r} -
A^{(2)}_{\tau_{N(t)}}|}{\sqrt{t}} ~~\rightarrow~~ 0~~~~~~{\rm
in}~~~~~{\rm probability}.
\eneq The arguments identical to the
ones employed in the proof of Lemma \ref{RetExpMom2} show that for
some positive $D$ and $\gamma$ and for any $n$
\beq \PROB \left\{
\sup_{0 \leq l \leq \tau_{n+1} - \tau_n }  |A^{(2)}_{\tau_n +l } -
A^{(2)}_{\tau_{n}} | > t \right\} \leq D e^{-\gamma t}~.
\eneq
Combined with (\ref{ab1}) this shows that for any $c>0$
\beq \label{ab2}
\sup_{s \in [\tau_{[N(t) -c \sqrt{t}]}, \tau_{[N(t) + c
\sqrt{t}]}] }
 \frac{|A^{(2)}_s -
A^{(2)}_{\tau_{N(t)}}|}{\sqrt{t}} ~~\rightarrow~~ 0~~~~~~{\rm
in}~~~~~{\rm probability}.
\eneq
By the CLT the quantity $ \frac{t - \tau_{N(t)}}{\sqrt{t}}$
converges to a normal random variable.
Therefore, by (\ref{ab2}) the second term in the right side of
(\ref{st}) tends to zero in probability. The first term tends to a
normal random variable by (\ref{ab0}). Q.E.D.

We shall extend Theorem \ref{CLT2pt} for a two--point
motion to an $n$-point motion after the next section.

\section{Exponential Mixing of the two--point process.}
\label{mixing}

In this section we prove exponential decay of correlations for
the two--point process, defined by stochastic flow (\ref{SF}),
starting outside of the diagonal $\Delta\subset M \times M$.
Recall  that by \cite{K} the two--point process has the unique
invariant measure $\mu^{(2)}$ on $M\times M \setminus \Delta. $
$\mu^{(2)}$ is a product measure, $\mu^{(2)}=\mu \times \mu$.

\begin{theorem}
\label{2ptMix} Let $B\in C^\infty (M\times M)$ be a function with
zero mean. For any point $z\in M\times M  \setminus \Delta $ let
$\rho_{z,B}(t)=E_z(B(z_t))$, where $z_t$ is the two--point motion
with $z_0=z=(x,y)$. Then for sufficiently small positive $\theta$
there are positive  $C(\theta)$ and $p(\theta)$, with the property
that $p(\theta)\rightarrow 0$ when $\theta \rightarrow 0$, such
that
\beq \label{corrdec}
|\rho_{z,B}(t)|\leq C e^{-\theta t}
\left(\frac{1}{d(x,y)}\right)^p.
\eneq
\end{theorem}

It turns out that it is sufficient to prove a similar statement
for $z \in \overline{G}_r$, that is for the initial point at least
distance $r$ away from the diagonal. \bprop \label{MixInd} With
notations of Theorem \ref{2ptMix}, for any small some positive $r$
there exist positive $\theta$ and $C$, such that if $u \in
\overline{G}_r$, then \beq |\rho_{u,B}(t)|\leq C e^{-\theta t}.
\eneq \eprop

Now, using Proposition \ref{MixInd}, we prove Theorem
\ref{2ptMix}.

{\it Proof of Theorem \ref{2ptMix}:}\ Let $\tau(z)$ be the
stopping time of hitting $ \overline{G}_r$, that is
$\tau(z) = \inf \{t: z_t \in \overline{G}_r \}$.
Then for the indicator function $\chi_{\{\tau <\frac{t}{2}\}}$ we
have
\beq \EXP_z B(z_t)=\EXP_z ( B(z_t)\chi_{\{\tau <\frac{t}{2}\}} )
+ \EXP_z ( B(z_t)\chi_{\{\tau  \geq \frac{t}{2}\}} ).
\eneq
To bound the first term, we use Markov property together with
Proposition \ref{MixInd},
\[
|\EXP_z ( B(z_t)\chi_{\{\tau <\frac{t}{2}\}} ) | \leq
\PROB_z \left\{ \tau <\frac{t}{2}\right\}
\sup_{s \in [\frac{t}{2}, t], u \in
\overline{G}_r} | \rho_{u,B}(s) | \leq
C \ \|B\|\ e^{-\frac{\theta t}{2}} ~.
\]
For
the second term we can apply Lemma \ref{Escape} and get that
\beq \nonumber
| \EXP_z ( B(z_t)\chi_{\{\tau  \geq \frac{t}{2}\}} )| \leq
\|B\|\ \PROB_z \left\{\tau(z) \geq \frac{t}{2}\right\} \leq C\
\|B\|\ \left(\frac{1}{d(x,y)}\right)^p e^{-\frac{\alpha t}{2}}~,
\eneq
where $p$ can be taken arbitrarily small by selecting a
sufficiently small $\alpha$. This
completes the proof of Theorem \ref{2ptMix}. Q.E.D.

{\it Proof of Proposition \ref{MixInd}:}\ Let $\hrho_{u,B}(\xi)$
be Fourier transform of $\rho_{u,B}.$
\beq
\beal \label{fourier}
\hrho_{u,B}(\xi)=\int_0^\infty e^{-i\xi t}\ \EXP_uB(u_t)\ dt =\
\EXP_u\left(\int_0^\infty e^{-i\xi t} B(u_t) dt\right).
\enal
\eneq
The main idea of the proof is to show that for each point
$u\in  \overline{G}_r$ the corresponding Fourier transform
function $\hrho_{u,B}(\xi)$ is analytic in $\xi$ in some strip
$|{\rm Im}\xi|<\alpha$, and $\int |\hrho_{u,B}(\sigma+is)|
d\sigma$ is uniformly bounded for $|s|<\alpha,~ u \in
\overline{G}_r$. Then by Paley-Wiener Theorem we have $\rho_{u,
B}(t) \leq c e^{-\theta t}$ for any $\theta < \alpha$ and  some
positive $c$. This clearly implies the Proposition.

Rewrite now the right-hand side of (\ref{fourier})
\beq
\label{fourier1} \EXP_u\left(\int_0^\infty e^{-i\xi t} B(u_t)
dt\right)= \EXP_u\left(\sum_{j=0}^\infty e^{-i\xi \tau_j}
\int_{\tau_j}^{\tau_{j+1}} e^{-i\xi(t-\tau_j)}B(u_t) dt\right),
\eneq
where $\tau_0 = 0$, and $\tau_j,~j\geq 1$ were defined in
Section \ref{pl}. Let
\beq \beta(u,\xi)=
\EXP_u\left(\int_{0}^{\tau_1} e^{-i\xi t} B(u_t) dt\right).
\eneq
Introduce on the space $C( \overline{G}_{r/2})$ of continuous
functions on $\overline{G}_{r/2}$ the transfer operator $\cR_\xi:
C( \overline{G}_{r/2}) \rightarrow   C( \overline{G}_{r/2})$
\beq
(\cR_\xi f)(u)= \EXP_u e^{-i\xi \tau_1} f(u_{\tau_1})~.
\eneq
By the Markov property
\[
\EXP_u\left(e^{-i\xi \tau_j} \int_{\tau_j}^{\tau_{j+1}}
e^{-i\xi(t-\tau_j)}B(u_t) dt\right) = [\cR_\xi^j
\beta(\cdot,\xi)](u)~.
\]
Now  using this formula and (\ref{fourier1}) we get \beq
\hrho_{u,B}(u,\xi)=\sum_{j=0}^{\infty} \left[\cR_{\xi}^{j}
\beta(\cdot,\xi)\right] (u)= \left[(1-\cR_{\xi})^{-1}
\beta(\cdot,\xi)\right] (u)~, \eneq provided that $\|\cR_{\xi}\| <
1$. Note the following properties of $\cR_{\xi} $ and
$\beta(\cdot,\xi)$

1) By Lemma \ref{RetExpMom}, $\cR_{\xi} $  and $\beta(\cdot,\xi)$
are  analytic in the half plane ${\rm Im} \xi < \gamma $. The
functions  $\beta(\cdot,\xi)$ and  $\beta'_{\xi} (\cdot,\xi)$ are
uniformly bounded in the same half plane.

2) Since $\tau_1 \in \{\delta, 2\delta, \dots \}$, the operator
$\cR_{\xi} $ is periodic with the period $2 \pi \delta^{-1}$. As
follows from the definition of $\cR_{\xi} $,
\[
\|\cR_{\xi} \| < 1, ~~~~~~~{\rm when} ~{\rm Im} \xi \leq  0, ~ \xi
\neq 2 \pi k \delta^{-1}~.
\]

3)  The operator $R_0$ is a generator of a mixing Markov chain.
Thus, in a neighborhood of $\xi = 2 \pi k \delta^{-1}$, the
operator $\cR_{\xi} $ can be written as
\[
\cR_{\xi}f = {\cS}_{\xi}  f + {\cT}_{\xi} f~.
\]
The operator $\cS_{\xi}$ is one dimensional,
  $\cS_{\xi}f = \lambda_{\xi} \mu_{\xi}(f) f_{\xi} $, where
$ \lambda_{\xi}$ is the largest (in absolute value) eigenvalue of
  $\cR_{\xi}$, while $ f_{\xi}$ and $ \mu_{\xi}$ are the
corresponding eigenfunctions of  $\cR_{\xi} $ and of the adjoint
operator, such that $ \mu_{\xi}(f_{\xi})=1$ . The operators
$\cS_{\xi}$ and $\cT_{\xi}$ commute, and $\|\cT_{\xi}\| < 1$. The
operators $cS_{\xi}$ and $cT_{\xi}$, and the vectors $f_{\xi}$ and
$\mu_{\xi}$ are analytic.

In a neighborhood of $\xi =  2 \pi k \delta^{-1}$ the operator
$(1-\cR_{\xi})^{-1}$ can be written as
\begin{equation}
\label{inverse}
(1-\cR_{\xi})^{-1} f = (1-\cT_{\xi})^{-1}f +
\frac{\cS_{\xi}f}{1-\lambda_{\xi}}~
\end{equation}
Thus, the operator $(1-\cR_{\xi})^{-1}$ admits a  meromorphic
continuation to  a half plane ${\rm Im} \xi < \alpha $ for some
$\alpha >0$, with the poles at the points $\xi =  2 \pi k
\delta^{-1}$. Note that $\lambda'_{\xi} \neq 0$ at the poles, and
therefore the poles are simple. From (\ref{inverse}) and Remark
(\ref{zer}) it follows that the residues of $(1-\cR_{\xi})^{-1}
\beta(\cdot,\xi)$ at $\xi =  2 \pi k \delta^{-1}$ are equal to
\[
\frac{\cS_{ 2 \pi k \delta^{-1}} \beta(\cdot,  2 \pi k
\delta^{-1})}{\lambda'_{\xi}}=
\frac{\EXP_{\mu_0}\left(\int_{\tau_0}^{\tau_1}
e^{-i  2 \pi k \delta^{-1} t} B(u_t) dt\right)\cdot
1}{\lambda'_{\xi}}=0~.
\]
Therefore $(1-\cR_{\xi})^{-1} \beta(\cdot,\xi)$ is holomorphic in
the half plane  ${\rm Im} \xi < \alpha $.

Next we show that
$|(1-\cR_{\xi})^{-1} \beta(\cdot,\xi)|$ is uniformly bounded in
the same half plane. Let $U$ be a neighborhood of $\xi = 0$, where
(\ref{inverse}) is valid. Since $\cR_{\xi}$ is periodic,
(\ref{inverse}) is valid in the same neighborhood of each of the
other poles of $(1-\cR_{\xi})^{-1}$. The function
$|(1-\cR_{\xi})^{-1} \beta(\cdot,\xi)|$ is bounded outside of the
union of such neighborhoods, since $(1-\cR_{\xi})^{-1}$ is a
bounded operator there. In the neighborhood $U_k$ of the point $ 2
\pi k \delta^{-1}$, the function $|(1-\cT_{\xi})^{-1}
\beta(\cdot,\xi)|$ is bounded, since  $(1-\cT_{\xi})^{-1}$ is a
bounded operator there. Finally, $|\frac{\cS_{\xi}\beta(\cdot,
\xi) }{1-\lambda_{\xi}}|$ is estimated by ${\rm Const}
|(\cS_{\xi}\beta(\cdot, \xi))'_{\xi}|$, which is bounded, since
$\cS_{\xi}$ is periodic, while  $\beta(\cdot, \xi)$ and
$\beta'_{\xi}(\cdot, \xi)$ are bounded.

\begin{lemma} \label{fourier3}
For any positive integer $N$ we have $\hrho_{u,B}(\xi)\xi^N\to 0 $
uniformly in $u \in G_{r/2}$,  as $\xi\to \infty, ~{\rm Im} \xi <
\alpha.$
\end{lemma}
{\it Proof:}\  We have $$|\hrho_{u,B}(\xi)\xi|=|(\partial_t
\rho_{u,B}(t)\hat)(\xi)|= |\hrho_{u,D(B)}(\xi)|$$ where $D$ is the
generator of the two point motion. Thus,  for any function $ B $
with zero mean, $|\hrho_{u, B}(\xi)| \leq \frac{K_1}{|\xi|}.$ In
particular $|\hrho_{u, D(B)}(\xi)|<\frac{K_1}{|\xi|}.$ Hence
$|\hrho_{u, B} (\xi)|  \leq \frac{K_2}{|\xi|^2}.$ Continuing by
induction we obtain the Lemma. Q.E.D.

Thus, $\hrho_{u,B}(\xi)$ is analytic in $|{\rm Im} (\xi)|<\alpha$
and $$\int |\hrho_{u,B}(\sigma+is)| ds$$ is uniformly bounded for
$|s|<\alpha$ by Lemma \ref{fourier3}. Thus Proposition
\ref{MixInd} is implied by Paley-Wiener Theorem. Q.E.D.

\section{An $n$-point motion} \label{multi}

Below we show the modifications needed to prove Theorem \ref{CLTnpt}.
We shall follow word-by-word the arguments given in case $n=2$ and
exhibit differences. For the proof of Theorem \ref{CLT2pt}
we use $\EXP_z(e^{\alpha(\tau_1-\sg_1)})<K$ or the fact that
time the two-point motion spends inside $r$-neighborhood of
the diagonal has exponential moments.

First we define an analog of return times of returning to
$r$-neighborhood of the diagonal $\Delta_r $ hitting
$r/2$-neighborhood $\Delta_{r/2} $ in between similarly to
(\ref{stop}). For the multipoint motion some care has to be taken.

The $r$-cylinder of the generalized diagonal and its outside
are given by
\beq \label{r-nb}
\beal
\Delta^{(n)}_r =\{(x^1,\dots,x^n)\in M\times \dots \times M:
\exists i<j\ d(x^i,x^j)=r\},\\
G^{(n)}_r =\{(x,\dots,x^n)\in M\times \dots \times M: \quad
\min_{i<j}d(x^i,x^j)> r\}.
\enal
\eneq
Denote $z^{(n)}_t=(x^1_t,\dots,x^n_t)$ the position of an $n$-point
motion at time $t$. Similarly to the two--point motion we decompose
the time interval $[0,t]$ into random time intervals as follows.
The initial  interval is $[0,\sg^{(n)}_1]$, where
$\sg^{(n)}_1=\sg^{(n)}_r(z^{(n)})$ is the first moment
$z^{(n)}_s \in \Delta^{(n)}_{r/2} $ the point $z^{(n)}_s$ hits
$\Delta^{(n)}_{r/2}$. Then we define two types of intervals of time.
For each $i<j$ formula (\ref{stop}) aplied to the two-point motion
$(x^i_t,x^j_t)$ defines the stopping times
$\{\sg^{ij}_s,\tau^{ij}_s\}_s$  of consequently hitting $\Delta_r$
and staying time $\dt$ inside $G_{r/2}$.
Consider the union of these intervals over all $i<j$
\beq
\cup_{i<j} \cup_{m\in \Z_+} [\sg^{ij}_m,\tau^{ij}_{m+1}].
\eneq
Renumerate this union as the one consisting of a countable number
of disjoint connected intervals. Denote renumerated intervals
by $\{[\sg^{(n)}_m,\tau^{(n)}_m]\}_{m\in \Z_+}$.  we call each
such an interval $[\sg^{(n)}_m,\tau^{(n)}_m]$ an interval of
the second type. During $s \in [\sg^{(n)}_m,\tau^{(n)}_m]$
the $n$-tuple $z^{(n)}_s$ stays
close to the generalized diagonal and the complement intervals
$\{[\tau^{(n)}_{m},\sg^{(n)}_{m+1}]\}_{m\in \Z_+}$ are of the
first type and correspond $z^{(n)}_s$ being in $G^{(n)}_{r/2}$
or $r/2$-away from the generalized diagonal. By compactness of
$G^{(n)}_{r/2}$ and nondegeneracy assumption (C$_n$) during
intervals of the first type the $n$-point motion is nondegenerate.

Using these stopping times we subdivide $[0,t]$ into
$[0,\sg^{(n)}_1],\ [\sg^{(n)}_1,\tau^{(n)}_1]$,
and so on $[\tau^{(n)}_{m(t)},t]$,
where this is the interval after the last
$\tau^{(n)}_{m(t)}$-return to $G^{(n)}_{r/2}$.

An additional difficulty arising for the $n$-point motion is
to prove an analog of Lemma \ref{Escape} for the return stopping
time $\tau^{(n)}_1$ defined above.

\begin{lemma} \label{REM} For flow (\ref{SF}) with the above
notations for the $n$-point motion
there exist positive $D_1$ and $\gamma_1$ such that for each
$u\in G^{(n)}_{r/2} $ and each $j\in \Z_+$ we have
\beq
\PROB\{\tau^{(n)}_1>t|\ z^{(n)}_0 = u \}
\leq D_1\ e^{-\gamma_1 t}.
\eneq
\end{lemma}
Suppose this Lemma is proven. Then implies that
$E_{z^{(n)}} (\exp \{\al \tau^{(n)}_1\})$ is finite for
some $\al>0$. Then the rest of the proof is the same proof as
for the $2$-point motion.

{\it Proof of Lemma \ref{REM}:}\ The proof is a little more
complicated than the one for Lemma \ref{RetExpMom}. The problem
is that $\tau^{(n)}_1$ the return time of the $n$-point motion
consists of two parts: first
$[0,\sg^{(n)}_1]$ when $z^{(n)}_{\sg^{(n)}_1}$ hits
$\Delta^{(n)}_{r/2} $ for the first time and
$[\sg^{(n)}_1,\tau^{(n)}_1]$ when $z^{(n)}_t$ is inside
$\Delta^{(n)}_r $, i.e. there is a pair of points
$x^i_t$ and $x^j_t$ which are inside $\Delta^{ij}_r$.
Exponential decay of the probability for
$\sg^{(n)}_1>t$ follows from compactness of $G^{(n)}_{r/2} $.
The problem which arises is that intervals of time
of return from  $\Delta^{ij}_{r/2} $ to $\Delta^{ij}_r $ for
different pairs $i$ and $j$ can be overlapping and, therefore,
Lemma \ref{REM} is not straightforward corollary from
Lemma \ref{Escape}.

The idea of the proof is to fix a pair of points $x^i$ and $x^j$ and
show that for a large time $t$ and a sufficiently small $r>0$
the probability that proportion of time $(x^i_s,x^j_s)$ spend in
the intervals of time of the first type $[\sg^{(n)}_k,\tau^{(n)}_k]$
compare to the whole time $t$ exceeds $1/n^2$ decays exponentially
in $t$, i.e. for some $K>0$ and $0<\theta<1$
\beq \label{expdecay}
\PROB_{(x^i,x^j)}
\left\{ \left| \sum_{\{r:\ \tau^{ij}_r(x^i,x^j)\leq t\}}
(\tau_r(x^i,x^j)-\sg^{ij}_r(x^i,x^j) \right| \leq
\frac{t}{n^2}\right\}\leq K \theta^t.
\eneq
Therefore, union over all pairs of point gives that for any
$z^{(n)}=(x^1,\dots, x^n)\in M\times \dots \times M
\setminus \Delta^{(n)} $ we have
\beq \nonumber
\beal
\PROB_{z^{(n)}}\left\{\tau_1^{(n)}(z^{(n)})-\sg_1^{(n)}(z^{(n)})>t
\right\}  \leq \\
\sum_{i<j} \PROB_{(x^i, x^j)}
\left\{\tau_1^{(2)}(x^i, x^j)-\sg_1^{(2)}(x^i, x^j)>
\frac{t}{n^2}\right\}.
\enal
\eneq
Therefore, if probabilities of every term on the right-hand side
is exponentially small in $t$, then probability in the left-hand
side is exponentially small. This would complete the proof of
Lemma \ref{REM}. Now we proof exponential decay of the right-hand
side probabilities of the two--point motion.

It is shown in \cite{BS} the two--point motion, defined by flow
(\ref{SF}), is ergodic and has a unique invariant measure $\mu^{(2)},$
such that $\mu^{(2)}(\Delta)=0$. Let $u \in G_{r/2}$. Denote by
$E_u(\tau_1(u))$ expected time of a point $u\in G_{r/2}$ hitting
$\Delta_{r/2}$, then $\Delta_r$, and then staying outside
$\Delta_{r/2}$ for a time at least $\dt$ for the two-point motion.
Then by ergodic theorem for stopping times $\tau_1(u)$ and $\sg_1(u)$,
defined above, we have
\beq \label{TimeRat}
\frac{\int E_u(\tau_1(u)-\sg_1(u)) dp(u)}{\int
E_u(\tau_1(u)) dp(u)}\to 1 \quad \textup{as} \quad
r\to 0,
\eneq
where $dp(u)$ is an invariant measure on $G_{r/2}$ induced by
flow (\ref{SF})
Let
\beq
\gamma(\xi_1,\xi_2,u,w)=\EXP_u(e^{\xi_1(\sg_2(u)-\tau_1(u))+
\xi_2(\tau_1(u)-\sg_1(u))}|z^{(0)}=u, z^{(1)}=w),
\eneq
with both $u, w \in G_{r/2}$. Define a tranfer operator
\beq
(\cM_\xi f)(u)=\int k(u,w) \gamma(\xi_1, \xi_2,u,w) f(w) dw,
\eneq
where $k(u,w)$ is the transition density of returns to $\Delta_r $.
Since $\cM_\xi$ depends analytically on $\xi=(\xi_1,\xi_2)$,
therefore, for small $\xi$ we have a decomposition
\beq
\cM_\xi f=\lambda_\xi \mu_\xi(f) f_\xi +
\Pi_\xi (f-\mu_\xi (f)f_\xi),
\eneq
where operators $\Pi_\xi$ and $\cL_\xi$ commute,
$\cL_\xi f_\xi=\lambda_\xi f_\xi$ and $\Pi_\xi$ is a contraction
operator, i.e. $\|\Pi_\xi\|<\theta$ for some $0<\theta<1$.
Direct calculation shows
\beq \nonumber
E_u\left(\exp\left[\xi\left(\sum_{j=1}^m n^2 [\tau_j(u)-\sg_j(u)]+
(\sg_{j+1}(u)-\tau_j(u)) \right)\right]\right)=
(\cM_\xi^m 1)(u)
\eneq
Since $\lb_0=1$ and using (\ref{TimeRat}) for any $n\in \Z_+$ one
can choose positive $r, K,$ and $\xi$ and some $0<\theta<1$ such
that for any positive integer $m$
\beq\nonumber
\EXP_u\left(\exp\left[\xi\left(\sum_{j=1}^m n^2 [\tau_j(u)-\sg_j(u)]+
(\sg_{j+1}(u)-\tau_j(u)) \right)\right]\right) \leq K \theta^m.
\eneq
By Chebyshev inequality this implies that
$$\PROB_u\left\{\frac{\sum_{j=1}^m [\tau_j(u)-\sg_j(u)]}{\tau_m}
\geq \frac{1}{n^2}\right\}\leq K \theta^m . $$
Similar argument shows that for any $T$ exceeding average return
time $\tau^*$ and any positive integer $m$ may be with a different
$K$ we get
$$\PROB_z(\tau_m>Tm)\leq K \theta^m. $$

Now we are ready to prove (\ref{expdecay}).
Let $m=[\frac{t}{K}].$
\beq
\beal
\PROB_{(x^i,x^j)} \left\{ \left|
\sum_{\{r:\ \tau^{ij}_r(x^i,x^j)\leq t\}}
(\tau^{ij}_r(x^i,x^j)-\sg^{ij}_r(x^i,x^j) \right|\leq\frac{t}{n^2}
\right\}
\leq \\
\PROB_{(x^i,x^j)}\{\tau_m>t\}+
\sum_{l=m+1}^\infty \PROB_{(x^i,x^j)}\left\{\frac{\sum_{j=1}^l
{(\tau_j-\sg_j)}}{\tau_l} \geq \frac{1}{n^2}\right\}.
\enal
\eneq
All the terms in the right-hand side decay exponentially with
exponent $\theta$. This completes the proof of Lemma \ref{REM}.
Q.E.D.

As we pointed out above the rest of the proof of Theorem \ref{CLTnpt}
is completely similar to the proof of Theorem \ref{CLT2pt} as we
pointed out above. Theorem \ref{CLTnpt} is proven. Q.E.D.

\section{Equidistribution.}
\label{curve}

In this section we  prove that images of smooth submanifolds
become uniformly distributed over $M$ as $t \to \infty$. More
generally, we prove that if  a measure $\nu$ has  finite
$p$-energy, defined in (\ref{energy}),  for some $p>0$, then
the image of this measure under the dynamics of stochastic flow
(\ref{SF}), satisfying conditions (A) through (E), becomes uniformly
distributed on $M$.

\begin{lemma} \label{CrEquid} Let $\nu$ be a measure on $M$ which
has finite $p$-energy for some $p>0.$ Let $b \in C^{\infty}(M)$
satisfy $\int b(x) d \mu (x) = 0$. Then there exist positive
$\gamma$ independent of $\nu$ and $b$, and $C$ independent of
$\nu$ such that for any positive $t_0$
\begin{equation} \label{res1}
\PROB \left\{ \sup_{t \geq t_0} \left|\int b(x_t) d\nu (x) \right|> C
I_p(\nu)^{1/2} e^{-\gamma t_0} \right\} \leq C  e^{-\gamma t_0}~.
\end{equation}
\end{lemma}

{\it Proof:}\ Without loss of generality we may assume that
$\nu(M) =1$ (otherwise we multiply $\nu$ by a constant).
 By the  exponential mixing of  the two--point
processes (Theorem \ref{2ptMix}) we have
\beq
\beal
\EXP\left(\int b(x_t) d\nu(x) \right)^2=\int\!\!\!\!\int_{M\times M}
\EXP_{(x,y)} (b(x_t) b(y_t)) d\nu(x) d\nu(y) \leq \\
C_1 \|b\|^2 e^{-\theta t} \int\!\!\!\!\int \frac{d\nu(x) d\nu(y)}{d^p(x,y)}
\leq C_2 I_p(\nu) e^{-\theta t}.
\enal
\eneq
Therefore,
\begin{equation} \label{mom1}
\EXP\left|\int b(x_t) d\nu(x) \right| \leq C_3 I_p(\nu)^{1/2}
e^{- \theta t/2} ~.
\end{equation}

This shows that the expectation of the integral decays
exponentially fast. Now we shall use standard arguments based on
Borel-Cantelli Lemma to show that $\int b(x_t) d \nu(x)$ itself
decays to zero exponentially fast.

Fix small positive $\alpha $ and $\kappa$, to be specified later.
Let $\tau_{n,j}=n+j e^{-\kappa n},$ $0 \leq j \leq e^{\kappa n}.$
By Chebyshev inequality, (\ref{mom1}) implies
\[
\PROB \left\{\left|\int b(x_{\tau_{n,j}} ) d\nu (x) \right| >
I_p(\nu)^{1/2} e^{-\alpha  n } \right\}
\leq C_3  e^{(\alpha- \theta /2) n }~.
\]
Taking the sum over $j$,
\begin{equation}
\label{ver1}
\PROB \left\{ \max_{0\leq j <e^{\kappa n}} \left|\int
b(x_{\tau_{n,j}} ) d\nu (x) \right| >I_p(\nu)^{1/2} e^{-\alpha  n }
\right\}
\leq C_3 e^{(\kappa + \alpha- \theta /2) n }~.
\end{equation}
Next we consider the oscillation of $\int b(x_t) d \nu(x) $ on the
interval $[\tau_{n,j}, \tau_{n,j+1} ]$. By Ito formula, due to
(\ref{SF}),
\begin{equation}
\label{ito1}
\beal \int b(x_t) d \nu(x) - \int b(x_{\tau_{n,j}}) d
\nu(x) = \\  \sum_{k =1}^d \int_{\tau_{n,j}}^t \int a_k(x_s) d \nu
(x) d\theta_k (s) + \int_{\tau_{n,j}}^t \int a_0(x_s) d \nu (x)
ds~,
\enal
\end{equation}
where $a_k$, $k=0,\dots,d$ are bounded functions on $M$. Each of
the integrals $\int_{\tau_{n,j}}^t \int a_k(x_s) d \nu (x)
d\theta_k(s)$ can be obtained from a Brownian motion by a random
time change. This time change has a bounded derivative, since
$a_k$ are bounded.  The absolute value of the  last term on the
right side of (\ref{ito1}) is not greater than $\frac{e^{-\kappa
n/4}}{2(d +1) }$ for sufficiently large $n$. Therefore, for
suitable $C_4, C_5, C_6> 0$, which are independent of $\nu$ since
$\nu(M) =1$,
\begin{equation}
\beal
\PROB \left\{ \sup_{\tau_{n,j} \leq t_1 \leq t_2  \leq
\tau_{n,j+1} } \left|\int b(x_{t_1} ) d\nu (x)  - \int b(x_{t_2} )
d\nu (x) \right| \geq e^{-\kappa n/4 } \right\} \leq \\ d \cdot
\PROB \left\{\sup_{0 \leq t_1 \leq t_2  \leq C_4 e^{-\kappa n}  }
|w(t_1) - w(t_2)| \geq \frac{ e^{-\kappa n/4 }}{d+1} \right\}
\leq
\\
C_5 \exp\left(-\frac{( \frac{e^{-\kappa n/4} }{d +1})^2 }{2 C_4
e^{-\kappa n} }\right) \exp(\frac{\kappa n}{2}) = C_5 \exp( -C_6
e^{\kappa n /2} ) \exp(\frac{\kappa n}{2})~.
\enal
\end{equation}
Combining this with (\ref{ver1}) we obtain
\[
\beal  \PROB \{  \sup_{n \leq t  \leq n+ 1 } \left|\int b(x_{t} )
d\nu(x) \right|
>I_p(\nu)^{1/2} (e^{-\alpha  n } + e^{-\kappa n /2}) \} \\
\leq C_3 e^{(\kappa + \alpha-\theta /2) n } + C_5 e^{\frac{3\kappa
n}{2}} \exp( -C_6 e^{\kappa n /2} )~. \enal
\]
This implies (\ref{res1}) if we take $\alpha = \kappa =\gamma =
\theta /10$ in the last inequality, and take the sum over all $n$
such that $n \geq t_0 - 1$. Q.E.D.

\section{CLT for measures}
\label{ScPathCLT}

\subsection{Energy estimate.}
\label{Aux}

 We prove the lemma which is needed in order to control the growth
of the the $p-$energy  of a measure. Let $\nu_t(\Omega)=\nu(x_t\in
\Omega).$ Also we shall write $\nu(f) = \int f d\nu.$
\begin{lemma}
\label{ExpEnergy} If $p$ is small enough, then for some $C>0$ $$
\EXP(I_p(\nu_t))\leq I_p(\nu)+C. $$
\end{lemma}
{\it Proof:}\ Fix a sufficiently small positive $r$. The goal is
to estimate from above
\beq
\EXP(I_p(\nu_t))=\EXP \ \int\!\!\!\!\int \frac{d\nu(x)\ d\nu(y)}{d^p(x_t,y_t)}.
\eneq
To get an upper
estimate we take care of when points $x_t$'s and $y_t$'s come
close one to the other. Consider two random sets:
\begin{eqnarray*}\nonumber
I_t^{0,-}&=&\{ (x,y)\in M \times M:\ \forall 0\leq s\leq t\ \ \
d(x_s,y_s)<r/2 \}, \\
I_t^{0,+}&=& M \times M \setminus I_t^{0,-}.
\end{eqnarray*}
For each pair of points $(x,y)$ from $I_t^{0,+}$ we define
the sequence of stopping times $\{ \sg_j(x,y)\}_{j\in \Z_+}$
(resp. $\{\tau_j(x,y)\}_{j\in \Z_+}$ with $\tau_j(x,y)>\sg_j(x,y)$
for all $j$) of hitting $r/2$-neighborhood (resp.
$r$-neighborhood) $\Delta_{r/2}$ (resp. $\Delta_{r}$) of the
diagonal. Now we subdivide the set $I_t^{0,+}$ into random subsets
\beq
\beal
I_t^{j,-}=\{ (x,y)\in M \times M:&\ \sg_j(x,y)\leq
t<\tau_j(x,y)\},
\\
I_t^{j,+}=\{ (x,y)\in M \times M:&\ \tau_j(x,y)
\leq t<\sg_{j+1}(x,y)\}.
\enal
\eneq
By Fubini's Theorem it suffices to
show that for some $C_1,\ C_2>0$ and any $t>0$ we have \beq \EXP\
d^{-p}(x_t,y_t) < C_1 d^{-p}(x,y) + C_2. \eneq Decompose  $\EXP
d^{-p}(x_t,y_t)$ into the sum
\beq
\EXP\
d^{-p}(x_t,y_t)=\sum_{j\in \Z_+} (\EXP\{ d^{-p}(x_t,y_t)\
\chi_{I_t^{j,+}}\}+ \EXP\{ d^{-p}(x_t,y_t)\  \chi_{I_t^{j,-}}\}),
\eneq
where $\chi_{I_t^{j,+}}$ and $\chi_{I_t^{j,+}}$ are the
characteristic functions. The first term is uniformly bounded
because $d(x_t,y_t)\geq r/2$. The second term means that at a time
$\sg_j(x,y)\leq t$ we have $d(x_{\sg_j(x,y)},y_{\sg_j(x,y)})=r/2$.
By renewal of the solution of the flow (\ref{SF}) we can apply
arguments from \cite{BS} (see proof of Theorem 3.19 on page 202) which
say that $d^{-p}(x_s,y_s)$ is a supermartingale, i.e.
$E\{ d^{-p}(x_t,y_t) \  \chi_{I_t^{j,-}}\}\leq (2/r)^p$.
This completes the proof.  Q.E.D.

\subsection{ Moment estimates and the Proof of the Main Result.}

This section is devoted to the proof of the main result of the
paper: CLT for the passive tracer (Theorem \ref{measure}). Recall
that we start with a nonrandom measure $\nu$ of finite $p$-energy
$I_p(\nu)<\infty$ for some $p>0,$ a stochastic flow of
diffeomorphisms (\ref{SF}), and an additive functional $\{A_t:\
t>0\}$ of the one--point motion, defined by the stochastic
differential equation (\ref{Add1}).

Let $\chi(t,\xi)=\nu(\exp\{\frac{i\xi}{\sqrt t} A_t\})$ be the
characteristic function of the functional $A_t(x)$. Below we shall
see that $\chi(t, \xi)$ is equicontinuous in $\xi$ when $\xi  \in
K, ~ t \in \N$ on a set of measure $1 - \delta$, where $\delta
> 0$ and a compact set $K$ are arbitrary. We shall further see
that
\begin{equation}
\label{fo2} \lim_{t \rightarrow \infty} \frac{\nu(|A_t -
A_{[t]}|)}{\sqrt{t}} = 0~~~~ {\rm almost ~~surely.}
\end{equation}
Finally we shall see that there exists $D(A)$, such that for any
$\xi$ fixed
\beq \label{fo3}
\lim_{n \rightarrow \infty} \chi(n,
\xi) = \exp\left( - \frac{D(A) \xi^2}{2}\right) ~~~{\rm
almost~~surely.}
\eneq
 Combining (\ref{fo3}) with (\ref{fo2})
above, and with the equicontinuity of $\chi$, we obtain Theorem
\ref{measure}. Thus it remains to verify the equicontinuity of
$\chi$ and formulas (\ref{fo2})-(\ref{fo3}).

In the next three lemmas  we estimate the growth rate of $A_t$ and
of its moments.
\begin{lemma} \label{L1}
Let $\delta, N_0$, and $N$ be arbitrary positive
numbers. There exists $C= C(\delta, N_0, N)$, such that for any
$t>0$, and any signed measure $\nu$, which satisfies $|\nu|(M)
\leq 1$, $I_p(|\nu|) \leq t^{N_0}$ we have
\[
\PROB\{ |\nu(A_t)| >  t^{\delta} \} \leq C t^{-N}~.
\]
\end{lemma}

{\it Proof:} Write the equation for $\nu(A_t)$ in Ito's form
\[
\nu(A_t) = \int_0^t \int \hat{a }(x_s) d\nu(x)ds + \sum_{k=1}^d
\int_0^t  \int \alpha_k(x_s) d\nu(x)  d\theta_k(s)~.
\]
By Lemma \ref{CrEquid}
\beq \label{sa1}
\PROB \{
\sup_{t^{\delta/3} \leq s \leq t} \left|\int b(x_s) d \nu(x) \right| >
e^{-\gamma t^{\delta/3}} t^{N_0/2} \} \leq C e^{- \gamma
t^{\delta/3}}~,
\eneq
where $b(x)$ is either $\hat{a }(x)$ or one
of $\alpha_k(x)$. The integrals
$\int_0^{t^{\delta/3}} \int \hat{a}(x_s) d \nu(x) d s $ and
$\sum_{k=1}^d \int_0^{t^{\delta/3}} [ \int \alpha_k(x_s) d\nu(x) ]
d\theta_k(s)$
are estimated using the facts that $\hat{a}$ and $\alpha_k$ are
bounded and that the stochastic integrals can be viewed as
Brownian motions with a random time change. The same integrals
over the interval $[t^{\delta/3}, t]$ are similarly estimated
using (\ref{sa1}). Q.E.D.

The proof of (\ref{fo2}) is similar to the proof of Lemma \ref{L1}
- one can write the expression for $A_t$ in Ito's form, and then
use the fact that the stochastic integral can be viewed as a time
changed Brownian motion. Alternatively, (\ref{fo2}) follows from a
more general result in \cite{LS}.

\begin{lemma}
\label{L2} There exists a constant $C > 0$, such that for any $k
\geq 0$ and any initial point $x$ \beq \label{es1} \PROB_{x}
\left\{ \frac{ | A_t |}{\sqrt{t}} >  k \right\} & \leq C
\exp\left(-\frac{k^2}{C}\right)~.~ \eneq
\end{lemma}
{\it Proof:}\ Recall that $A_t = \int_0^t  \hat{a }(x_s) ds +
\sum_{k=1}^d \int_0^t  \alpha_k(x_s)  d\theta_k(s)$. For the first
term by the large deviation theory
\beq \nonumber
\PROB_{x} \left\{ \frac{|\int_0^t \hat{a} (x_s) ds|}{\sqrt{t}} > k
\right\}
\leq C \exp\left(-\frac{k^2}{C}\right)\ \ \textup{for some}\ \
C>0,
\eneq
since $x_s$ is a mixing diffusion process on a compact manifold.
For each of the stochastic integrals, recall that $\int_0^t
\alpha_k(x_s) d\theta_k(s)$ can be viewed as a time changed
Brownian motion, with the derivative of the time change bounded.
Therefore,
\beq
\beal
\PROB_{x} \left\{ \frac{|\int_0^t \alpha_k (x_s) d
\theta_k(s)|}{\sqrt{t}} > k \right\} & \leq
\\
\PROB \left\{ \frac{\sup_{s \leq c t} | W_s|}{\sqrt{t}} > k \right\}
\leq  C & \ \exp\left(-\frac{k^2}{C}\right)\ \ \textup{for some}\ \
C>0.
\enal
\eneq
Therefore the estimate (\ref{es1}) holds, with possibly a
different constant $C$.
 Q.E.D.

\begin{lemma}
\label{L3} For any positive $\delta$ and any $ N \in \N$ there
exists a constant $C > 0$ such that for any measure $\nu$ with
$|\nu|(M) \leq 1$  and any $n \in \N$
\beq \label{es2}
\PROB \left\{ |\nu|(| A_t|^n)
> n!  t^{(\frac{1}{2} + \delta)n } \right\} \leq
C t^{-N-\frac{\delta n}{2}}  ~.~
\eneq
\end{lemma}
{\it Proof:}\ Without loss of generality we may assume that $\nu$
is a probability measure.  By Jensen's inequality $(\int|A_t|^n d
\nu)^l \leq \int |A_t|^{nl} d\nu$ for $l \in \N$. It is therefore
sufficient to estimate the probability $\PROB \{ \nu(|A|^{nl}) >
(n!)^l t^{(\frac{1}{2} + \delta)n l} \}.$ By Chebyshev inequality
\beq \label{ch}
\PROB \left\{ \nu\left(|A|^{nl}\right) > (n!)^l t^{(\frac{1}{2} +
\delta)n l} \right\} \leq \frac{\int \EXP_{x} |A|^{nl} d \nu(x) }{(n!)^l
t^{(\frac{1}{2} + \delta ) n l}}~.
\eneq
Take $l > \frac{2N}{\delta}$. Then the right side of (\ref{ch})
is not greater than $\frac{\sup_x \EXP_{x} |A|^{nl} }{(n!)^l
t^{\frac{1}{2} n l}} t^{-N - \frac{\delta n}{2}}.$ This is less
than $C t^{-N - \frac{\delta n}{2}}$ by Lemma \ref{L2}.
 Q.E.D.

Put $n_t=[t^{1/3}]$, $\tau_t=t/n_t$, and for each $0<s<t$ denote
the increment of the functional $A_t(x)$ from time $s$ to time $t$
by \beq \Delta_{s,t}(x)=A_t(x)-A_s(x). \eneq We split the time
interval $[0,t]$ into $n_t$ equal parts and decompose \beq
\label{part} A_t(x)=\sum_{j=0}^{n_t-1}
\Delta_{j\tau_t,(j+1)\tau_t}(x). \eneq The idea is to prove that
this is a sum of weakly dependent random variables and that the
CLT holds for almost every realization of the Brownian motion. We
need an estimate on the correlation between the inputs from
different time intervals. For any positive $\tau, s,$ and $l$ with
$\tau \leq l$ we denote
\[
\nu(\Delta_{l-\tau,l} \Delta_{l,l+s}) = \int
\Delta_{l-\tau,l}(y)\Delta_{l,l+s}(y)d \nu(y).
\]
\begin{lemma}
\label{L4} Let some positive $c_1, c_2, \gamma_1,$ and $\gamma_2$
be fixed, and consider $s$ and $\tau$ which satisfy $c_1
l^{\gamma_1} \leq s, \tau \leq c_2 l^{\gamma_2}$. For any positive
$\delta$,  $N$, and $N_0$ there exists a constant $C$ such that
for $l \geq C$, and any measure $\nu$, which satisfies $|\nu|(M)
\leq 1$, $I_p(|\nu|) \leq l^{N_0}$ we have \beq \PROB \{
|\nu(\Delta_{l-\tau, l} \Delta_{l, l+s})|
> l^{\delta} \} \leq l^{-N}~. \eneq
\end{lemma}
{\it Proof:} Without loss of generality we may assume that $\nu$
is a probability measure.  We start by decomposing each of the
segments $[l-\tau, l], [l, l+s]$ into two:
\[
[l-\tau, l] = \Delta_1 \cup \Delta_2 =  [l-\tau, l- \ln^2 l] \cup
[l -\ln^2 l, l]~,
\]
\[
[l, l+s] =  \Delta_3 \cup \Delta_4 = [l, l+  \ln^2 l] \cup [l +
\ln^2 l, l+s]~,
\]
We denote
\[
\bar{\Delta}_{a,b} = \Delta_{a,b} \chi_{\{\Delta_{a,b} \leq
(b-a)^2 \} }~.
\]
By Lemma \ref{L2} there is a constant $C$  such that
\[
(\EXP_{x} \Delta_{a, b}^2)^{1/2} \leq C (1+b-a)~~~~{\rm and}
\]
\[
(\EXP_{x} (\Delta_{a, b} - \bar{\Delta}_{a,b})^2)^{1/2} \leq C
e^{-\frac{b-a}{C}}~.
\]
Therefore for two segments $[a,b]$ and $[c,d]$ such that $b \leq
c$
\[
|\EXP_{x}(\Delta_{a,b} \Delta_{c,d} - \bar{\Delta}_{a,b}
\bar{\Delta}_{c,d})| = |\EXP_{x}((\Delta_{a,b}  -
\bar{\Delta}_{a,b}) {\Delta}_{c,d}) +
\EXP_{x}(\bar{\Delta}_{a,b}(\Delta_{c,d} - \bar{\Delta}_{c,d}))|
\leq
\]
\[
C(1+|b-a| + |d-c| )(e^{-\frac{b-a}{C}} + e^{- \frac{d-c}{C}})~.
\]
After  using Chebyshev's inequality we obtain that for any
positive $k$ \beq \label{u1} \beal \PROB \{ |\nu(\Delta_{a,b}
\Delta_{c,d}) - \nu (\bar{\Delta}_{a,b} \bar{\Delta}_{c,d})| >k \}
\leq
\\
\frac{C(1+ |b-a| + |d-c|)(e^{-\frac{b-a}{C}}+ e^{-
\frac{d-c}{C}})}{k} \enal \eneq In the same way one obtains \beq
\label{u2} \beal \PROB \{ |\nu(\Delta_{a,b} \Delta_{c,d}) - \nu
(\bar{\Delta}_{a,b} {\Delta}_{c,d})| >k \} \leq
\\
\frac{C(1+|b-a|+|d-c|)(e^{-\frac{b-a}{C}}+e^{- \frac{d-c}{C}})}{k}
\enal \eneq The contribution from $\nu(\Delta_2 \Delta_3)$ is
estimated using estimate (\ref{u1}) with $k = l^{\delta}$ : for
any positive $ \delta$ and $N$ for sufficiently large $l$
\[
\PROB \{ |\nu(\Delta_2 \Delta_3)| >  l^{\delta} \} \leq l^{-N}~.
\]
The contribution from each of the other three products is
estimated using the fact that the segments are separated by a
distance of order $\ln^2 l$. Let us for example prove that $\PROB
\{ |\nu(\Delta_1 \Delta_3)| >  l^{\delta} \} \leq l^{-N}.$ By
taking $k = l^\delta$ in (\ref{u2}) we obtain
\[
\PROB \{ |\nu(\Delta_1 \Delta_3)-\nu(\bar{\Delta}_1 \Delta_3)| >
l^{\delta} \} \leq l^{-N}~.
\]
In order to estimate $\nu(\bar{\Delta}_1 \Delta_3)$ we apply the
change of measure \beq \label{chm} \nu(\bar{\Delta}_1 \Delta_3) =
\int \bar{\Delta}_{l-\tau, l - \ln^2 l} \Delta_{l, l+ \ln^2 l} d
\nu(x) = \int \Delta_{ \ln^2 l, 2 \ln^2 l} (x) d \hat{\nu}(x)~,
\eneq where
\[
\hat{\nu}(A) = \int \chi_{ \{ x_{l - \ln^2 l} \in A \} }
\bar{\Delta}_{l-\tau, l - \ln^2 l} (x) d\nu(x)~.
\]
The measure $\hat{\nu}$ has a density with respect to $\nu_{l -
\ln^2 l}$ which is bounded by $l^2$ since $\bar{\Delta}_{l-\tau, l
- \ln^2 l} $ is bounded. Therefore by Lemma \ref{ExpEnergy} for
any positive $N$ there exists $M$ such that
\[
\PROB \{I_p (\hat{\nu}) > l^M \} \leq l^{-N}~.
\]
The right side of (\ref{chm}) is written as
\beq \nonumber
\beal
\int \Delta_{ \ln^2 l, 2 \ln^2 l} (x) d \hat{\nu}(x)= \int_{\ln^2
l}^{2 \ln^2 l}  \int \ha(x_s) d \hat{\nu}(x) ds+
\\
\sum_{k=1}^d \int_{\ln^2 l}^{2 \ln^2 l} \int \alpha_k(x_s)
d \hat{\nu}(x) d\theta_k(s)~.
\enal
\eneq
We now proceed as in the proof of Lemma \ref{L1}. The contribution
from $\nu(\Delta_1 \Delta_4)$ and $\nu(\Delta_2 \Delta_4)$ is
estimated in exactly the same way. Q.E.D.

Our next statement concerns the asymptotic behavior of the second
moment of the functional $A_t$
\begin{lemma} \label{L5}
The following limit exists and the convergence is uniform in the
initial point $x$
\beq
D(A) = \lim_{n \rightarrow \infty}
\frac{\EXP_{x} (A_n^2)}{n}~.
\eneq
\end{lemma}
{\it Proof:} Let $\gamma_j = A_{j+1} - A_j$. Then \beq \label{m1}
\quad \quad \frac{\EXP_{x}(A_n)}{n} = \frac{ \sum_{i = 0}^n
\EXP_{x}(\gamma_i^2)}{n}+ \frac{ 2 \sum_{i=0}^n \sum_{r = 1}^{n-i}
\EXP_{x}(\gamma_i \gamma_{i+r})}{n}. \eneq The following two
statements imply the existence of the limits for the two sums on
the right hand side of (\ref{m1})

(a) $|\EXP_{x}(\gamma_i \gamma_{i+r})| \leq C e^{-\delta r}~~$ for
some positive $C$ and $\delta$.\

(b) $\lim_{i \rightarrow \infty} \EXP_{x}(\gamma_i \gamma_{i+r}) =
f_r~~$ uniformly in $x$ for some $f_r$.

To prove (a) we write \beq \nonumber \beal |\EXP_{x}(\gamma_i
\gamma_{i+r})| = |\EXP_{x}(\gamma_i E(\gamma_{i+r} | {\mathcal
F}_{i+1}))| \\ \leq \left\{\EXP_{x}(\gamma_i^2)\right\}^{1/2}
\left\{\sup_y |\EXP_y(\gamma_{r-1} )|\right\}~. \enal \eneq Here
we assumed that $r \geq 1$. The case $r=1$ follows from the next
estimate.  The first factor is estimated as
\[
\left\{\EXP_{x}(\gamma_i^2)\right\}^{1/2} \leq \sup_y
\left\{\EXP_y(A_1^2)\right\}^{1/2} = C~.
\]
The last expectation is equal to
\[
|\EXP_y(\gamma_{r-1} )|= | \EXP_y \int_{r-1}^{r} \hat{a}(y_t) dt |
\leq c e^{-\delta r}~,
\]
where the inequality is due to the exponential mixing of the
one--point process. This proves (a).

Let $\phi_r(x) = \EXP_x(\gamma_0 \gamma_{r})$. The following limit
exists and is uniform in $x$ by the ergodic theorem \beq
\label{expectation} \beal \lim_{i \rightarrow \infty} \EXP_{x}
(\gamma_i \gamma_{i+r})= \lim_{i \rightarrow \infty} \EXP_{x}(\EXP
\gamma_i  \gamma_{i+r} |{\mathcal F}_{i})
\\
=
\lim_{i \rightarrow \infty} \EXP_{x} (\phi_r(x_i))= \int \phi_r(y)
d \mu(y)=f_r~. \enal \eneq This proves (b). Q.E.D.

The next lemma provides a linear bound (in probability) on the
growth of $\nu(A_t^2)$. Note that such a bound implies the
equicontinuity of $\chi(t,\xi)$ in the sense discussed above.

\blm \label{L6} For any positive $N, N_0$, and $\rho>0$ there is
$C>0$ such that for any measure $\nu$ which satisfies $|\nu|(M)
\leq 1$, $I_p(|\nu|) \leq t^{N_0}$ we have
 \beq \nonumber \beal \PROB
\left\{ \left| \nu(A_t^2) - \nu(M) D(A) t \right| \geq \rho t
\right\} \leq C\ t^{-N}. \enal \eneq \elm

{\it Proof:}\ Let us prove that  $ \PROB \{  \nu(A_t^2) - \nu(M)
D(A) t  \geq \rho t \} \leq Ct^{-N}$.
 The estimate
from below can be proved similarly. Consider the event $Q_t$ that
$I_p(\nu_{j\tau_t}) \leq t^{N_1}$ for all $j < n_t$. By Lemma
\ref{ExpEnergy}
we have $\PROB\{Q_t\} \geq 1 - t^{-N}$ if a
sufficiently large $N_1$ is selected. Fix any $\delta$ with $0 <
\delta < 1/20$.
 Let $R_t$ be the event that
 $\nu(\Delta_{j\tau_t,(j+1)\tau_t}^2)\leq\tau_t^{1+\delta}$ for all
 $j < n_t$. By  Lemma \ref{L3} we have $\PROB\{R_t\} \geq 1 -
t^{-N}$. Let $\beta_j= \nu( \Delta_{(j-1)\tau_t, j\tau_t}^2)
\chi_{\{Q_t \cap R_t \}}$ and $\mathcal B_j=\sum_{k=1}^j \beta_k.$
We shall prove that
 \beq  \label{h1} \EXP\left\{\exp \left(
\frac{\mathcal B_j -j(\nu(M) D(A)+\rho)\tau_t}{t^{5/6}}
\right)\right\}\leq (1-\rho t^{-1/6}/2)^j~.
 \eneq
The proof will proceed by induction on $j$. First we show that
 \beq  \label{h2} \EXP\left\{\exp \left(
\frac{\mathcal \beta_j -(\nu(M) D(A)+\rho)\tau_t}{t^{5/6}}
\right)\right\}\leq (1-\rho t^{-1/6}/2)~.
 \eneq
Indeed, using Taylor expansion
\beq
\beal
\EXP\left\{\exp \left(
\frac{\mathcal \beta_j -(\nu(M) D(A)+\rho)\tau_t}{t^{5/6}}
\right)\right\} = 1 +
\\
\EXP  \frac{ \beta_j -(\nu(M)
D(A)+\rho)\tau_t}{t^{5/6}} +
\EXP \sum_{k = 2}^{\infty} \left(
\frac{\mathcal \beta_j -
(\nu(M) D(A)+\rho)\tau_t}{t^{5/6}}\right)^k/k!~.
\enal
\eneq
 Since by definition
$\beta_j\leq \nu(\Delta_{j\tau_t,(j+1)\tau_t}^2)$, we have $\EXP
\beta_j \leq \EXP \nu( \Delta_{j\tau_t,(j+1) \tau_t}^2).$ From
Lemma \ref{L6} it easily follows that $\EXP \nu(
\Delta_{j\tau_t,(j+1) \tau_t}^2) \leq (\nu(M) D(A)+\rho/4)\tau_t$
for large $t$. The expectation of the infinite sum is less than
$\rho t^{-1/6}/4$ since
 $\beta_j \leq \tau_t^{1+\delta}$. This proves
(\ref{h2}).

Assume that (\ref{h1}) holds for some $j$. Then,
\beq \label{h3}
\beal
\EXP & \left\{ \exp\left(\frac{\mathcal
B_{j+1}-(\nu(M)D(A)+\rho) (j+1) \tau_t}{t^{5/6}}\right)\right\} =
\\
\ & \EXP\left\{ \exp\left(\frac{\mathcal B_{j}-(\nu(M)D(A)+\rho) j
\tau_t}{t^{5/6}}\right) \right.
\\
 \EXP & \left.\left( \exp\left(
\frac{\beta_{j+1}-(\nu(M)D(A)+\rho) \tau_t}{t^{5/6}}
\right)\ |\ {\mathcal F}_{j\tau_j}\right)\right\}~.
\enal
\eneq
Due
to (\ref{h2}) and since $I_p(\nu_{j\tau_t}) \leq t^{N_1}$ on $Q_t$
by the by the Markov property the
 conditional expectation on the right side of (\ref{h3}) is not greater
than $1-\rho t^{-1/6}/2$. Therefore, \beq \label{h4} \beal
\EXP\left\{ \exp\left(\frac{\mathcal B_{j+1}-(\nu(M)D(A)+\rho)
(j+1) \tau_t}{t^{5/6}}\right)\right\} \leq
\\
(1-\rho t^{-1/6}/2) \EXP\left\{ \exp\left(\frac{\mathcal
B_{j}-(\nu(M)D(A)+\rho) j \tau_t}{t^{5/6}}\right)\right\}~. \enal
\eneq This proves (\ref{h1}). It follows from (\ref{h1}) with $j =
n_t$ that

\[
\PROB \{  \mathcal B_{n_t} - \nu(M) D(A) t  \geq \rho t \} \leq
Ct^{-N}
\]
Recall that
\[
\PROB \left\{
 \mathcal  B_{n_t} \neq  \sum_{j=0}^{n_t-1}
\nu(\Delta_{j\tau_t,(j+1)\tau_t}^2) \right\}\leq C t^{-N} .
\]
by Lemmas \ref{ExpEnergy} and \ref{L3}.

Finally, direct application of Lemma \ref{L5} to pair products of
$\Delta$'s at different time segments gives that for any positive
$\delta$ we have
 \beq  \PROB \left\{\left| \nu(A^2_t)- \sum_{j=0}^{n_t-1}
\nu(\Delta_{j\tau_t,(j+1)\tau_t}^2)\right|  > t^{1/3+\delta} \right\} \leq C
t^{-N}  \eneq for sufficiently large $t$. This completes the proof
of the lemma. Q.E.D.

\begin{lemma} \label{charact}  Let $\rho,  N_0$, and $N$ be arbitrary
positive numbers. There exists a constant $C>0$, such that for any
$t>0$, and any signed measure $\nu$, which satisfies $|\nu|(M)
\leq 1$, $I_p(|\nu|) \leq t^{N_0}$ we have \[ \PROB  \left\{
\left| \nu \left(\exp\left\{ \frac{i\xi}{\sqrt t}
\Delta_{0,\tau_t} \right\} \right) - \nu(M) \left(1-\frac{D(A)
\xi^2}{2 t^{1/3}}\right) \right| \geq \rho  t^{-1/3}\right\} \leq
C t^{-N}~.
\]
\end{lemma}

{\it Proof Lemma \ref{charact}:}\ Consider the Taylor expansion of
the \\ function $\exp(\frac{i\xi}{\sqrt{t}}\Delta_{0,\tau_t}(x))$
\beq  \beal \nu \left( \exp\left\{\frac{i\xi}{\sqrt{t}}
\Delta_{0,\tau_t}\right\}\right) = \\ \nu(M)
+\frac{i\xi}{\sqrt{t}} \nu \left( \Delta_{0,\tau_t}\right) - &
\frac{\xi^2}{2t} \nu \left(\Delta^2_{0,\tau_t}\right)+
\sum_{k=3}^{\infty}(\frac{i\xi}{\sqrt{t}})^k
 \nu \left(\Delta^k_{0,\tau_t}\right)~ \enal \eneq
 By Lemma \ref{L1}  for any
$\delta >0$ we have
 \beq
\PROB \left\{ \left| \frac{i\xi}{\sqrt{t}} \nu \left(
\Delta_{0,\tau_t}(x)\right) \right| > t^{-1/2+\delta} \right\} < Ct^{-N}~.
\eneq
 By Lemma
\ref{L6} almost certainly we have
\beq \left| \nu
\left(\Delta^2_{0,\tau_t}(x)\right) -\nu(M) D(A)\tau_t\right|\leq
\rho \tau_t.
\eneq
To estimate the tail we apply Lemma \ref{L3}.
This proves the lemma. Q.E.D.

{\it Proof Theorem  \ref{measure}:} It remains to demonstrate
that (\ref{fo3}) holds.  First we show that
\beq  \label{sta}
\beal
 \PROB  \left\{ \left| \nu \left( \exp \{ \frac{i\xi}{\sqrt t}
\Delta_{0,(j+1)\tau_t}\} \right) - \left( 1-\frac{D(A)\xi^2}{2
t^{1/3}} \right) \nu \left( \exp \{ \frac{i\xi}{\sqrt t}
\Delta_{0,j\tau_t} \} \right) \right|  \right.\\ \left.
\geq \rho t^{-1/3} \right\} \leq  C t^{-N}.
\enal
\eneq
Write
\beq
\beal
\label{sta1}
 \nu \left( \exp \left\{ \frac{i\xi}{\sqrt t}
\Delta_{0,(j+1)\tau_t}\right\} \right)=  \nu \left( \exp \left\{
\frac{i\xi}{\sqrt t} \Delta_{0,j\tau_t}\right\}
 \exp \left\{
\frac{i\xi}{\sqrt t} \Delta_{j\tau_t,(j+1)\tau_t }\right\}\right)=\\
 \hat{\nu} \left( \exp \left\{ \frac{i\xi}{\sqrt t}
\Delta_{0,\tau_t}\right\} \right)~,
\enal
 \eneq
where $\hat{\nu}$ is a random measure, defined by
\[
\hat{\nu}(A) = \int  \exp \left\{ \frac{i\xi}{\sqrt t} \Delta_{0,j
\tau_t}\right\} \chi_{ \{x_{j \tau_t} \in A \} } d \nu(x)~.
\]
Note that by Lemma \ref{ExpEnergy} for some $N_0$
\[
\PROB \{ I_p(|\hat{\nu}|) > t^{N_0} \} \leq C t^{-N}~.
\]
Thus the right side of (\ref{sta1}) can be estimated with the help
of Lemma \ref{charact}
 \[ \PROB  \left\{
\left| \hat{\nu} \left(\exp\left\{ \frac{i\xi}{\sqrt t}
\Delta_{0,\tau_t} \right\} \right) - \hat{\nu}(M)
\left(1-\frac{D(A) \xi^2}{2 t^{1/3}}\right) \right| \geq \rho
t^{-1/3}\right\} \leq C t^{-N}~.
\]
This is exactly the same as (\ref{sta}). Applying (\ref{sta})
recursively for $j=n_t-1,...,1$ we obtain that for any positive
$N$ and $\rho$ there is $C > 0$ such that
\beq \nonumber
\beal
\PROB \{ (1-\frac{D(A) \xi^2 + \rho}{2t^{1/3}})^{n_t} \leq
 \nu \left( \exp \{ \frac{i\xi}{\sqrt t}
\Delta_{0,j\tau_t}\} \right)
\\
\leq (1-\frac{D(A) \xi^2 -\rho}{2t^{1/3}})^{n_t} \}
\geq 1 - Ct^{-N}~.
\enal
\eneq
Therefore,
\[
\PROB \{\exp (-\frac{D(A) \xi^2 + \rho}{2}) \leq
 \chi(t, \xi) \leq
 \exp (-\frac{D(A) \xi^2 - \rho}{2})  \} \geq 1 - Ct^{-N}~.
\]
This implies (\ref{fo3}), which completes the proof of the
theorem. Q.E.D.

\section{The Dissipative Case}
\label{ScDiss}

In this section we extend our CLT for measures
(Theorem \ref{measure}), proved in the previous section
for measure-preserving stochastic flows (see condition (A))
to the dissipative case. In other words, we consider
stochastic flows defined by the stochastic differential
equation (\ref{SF}) satisfying conditions (B), (C$_2$), (D),
and (E) the sum of whose Lyapunov exponents is negative.
First, notice that without measure-preservation assumption
it is no longer true that generically  the largest Lyapunov
exponent is positive. However, the case when all the Lyapunov
exponents are negative is well understood \cite{L3}, so
we shall concentrate on the case with at least one positive
eigenvalue. The main result of this section is CLT for measures
(Theorem \ref{ThDCLT} below).

Let $m$ be the invariant measure of the one--point process,
which is unique by hypoellipticity assumption (B).
Let $m_2$ be the invariant measure of the two-point process
which is supported away from the diagonal. Such a measure
exists and is unique for the processes with positive largest
exponent by the results of \cite{BS}. Moreover we have exponential
convergence to $m_2$ since the proof of Theorem \ref{2ptMix} never used the
assumption of volume preservation.

\begin{theorem}\label{dissip}
With the notations above there is a family of probability measures
$\{\mu_t:\ t>0\}$ such that
\newline
(A) For any measure $\nu$ of finite $p$-energy
$$ \lim_{n\to\infty} \phi_{-n,t}^* (\nu)=\mu_t $$
almost surely.
\newline
(B)\ the process $t \to \mu_t$ is Markovian and push forward $\phi_t$
by the time $t$ stochastic flow (\ref{SF}) satisfies
$\phi_t^*(\mu_0)=\mu_t$;
\newline
(C) For any continuous function $A$ for any measure $\nu$ of positive
$p$--energy
$$ \left| \nu(A(x_t))-\mu_t(A) \right|\to 0 \quad
as \quad t\to+\infty$$
almost surely.
\end{theorem}
\begin{remark}
The measures $\mu_t$ and especially their dimensional
characteristics were studied in several papers \cite{LY1, LY2, L1,
L2, L3}. The questions which we discuss here are different from
the  ones studied in these papers and we will not use any of their
results.
\end{remark}

Let $\{A_t:\ t>0\}$ be an additive functional of the one--point motion
satisfying (\ref{Add1}) and (\ref{CenterOfMass}).
Denote by $\bar a(t)=\mu_t(a)$ and $\bar\al_k(t)=\mu_t(\al_k)$
for $k=1,\dots,d$ averages of $a$ and $\al_k$'s with respect to
$\mu_t$ and define two additive functionals
\beq \label{decomposition}
\beal
dC_t= & \sum_{k=1}^d \bar \al_k(t) \circ d\theta_k(t)+
\bar a(t)dt, \\
dB_t(x)= & \sum_{k=1}^d (\al_k(x_t)-\bar \al_k(t)) \circ
d\theta_k(t)+ (a(x_t)-\bar a(t))dt. \enal \eneq Note that by the
standard theory of Markov processes $C_t^*=\frac{C_t}{\sqrt{t}}$
is asymptotically Gaussian with zero mean and some variance
$D''(A).$\footnote{It is easy to see from Theorem \ref{dissip}(A)
that $\forall a\in C(M)$ $\forall t$ we have $\EXP(\mu_t(a))=m(a)$
so the vanishing of the mean follows from (\ref{CenterOfMass}).}

\begin{theorem}
\label{ThDCLT} Let $\cM^{\theta,*}_t$ be the measure on $\R$
defined on Borel sets $\Omega \subset \R$ by \beq \label{Bmera}
\cM^{\theta,*}_t(\Omega)=\nu\left\{x\in M:\
\frac{B^\theta_t(x)}{\sqrt t}\in\Omega\right\}. \eneq Then there
is a constant $D'(A)$ such that almost surely $\cM^{\theta,*}_t$
converges weakly to a Gaussian measure with zero mean and variance
$D'(A).$
\end{theorem}
We can reformulate Theorem \ref{ThDCLT} as follows.
\begin{corollary}
\label{dis-measure}
(a) Almost surely for large $t$ the measure $\cM_t^\theta$
defined by (\ref{indmeas})
is asymptotically Gaussian with a random drift $C_t^*$ and deterministic
variance $D'(A).$
\newline
(b) As $t\to+\infty$ the distribution of the drift is asymptotically
Gaussian with zero mean and the variance $D''(A).$
\end{corollary}

{\it Proof of Theorem \ref{dissip}:}
Given a measure $\nu$ of finite $p$-energy denote
$$\mu_t^{(n)}(\nu)=\phi_{-n,t} \nu.$$
\begin{lemma}
\label{LmBackMix}
There is a constant $\rho<1$ such that
$\forall A\in C(M)$ $\forall t$ $\forall \nu_1, \nu_2$ of
finite $p$ energy almost surely there exists a constant $C=C(\theta)$
such that
\begin{equation}
\label{BackMix}
\left| \mu_t^{(n)}(\nu_1)(A)-\mu_t^{(n)}(\nu_2)(A)\right|\leq C \rho^n.
\end{equation}
\end{lemma}
\proof
\beq
\beal
\left| \EXP\left(
\mu_t^{(n)}(\nu_1)(A)-\mu_t^{(n)}(\nu_2)(A)\right)\right|= \\
\left| \EXP(A(x_t)) d\nu_1(x_{-n})-
\EXP(A(x_t)) d\nu_2(x_{-n})\right|
\leq  \Const\rho_1^n
\enal
\eneq
since both terms are exponentially close to $m(A)$ by the exponential
mixing of
one--point process. Likewise
\beq \label{BackVar}
\beal
\EXP\left(\left[\mu_t^{(n)}(\nu_1)(A)-\mu_t^{(n)}(\nu_2)(A)
\right]^2\right)= \\
\int\!\!\!\!\int A(x_t) A(y_t) d(\nu_1\times\nu_1)(x_{-n}, y_{-n})+
\\
\int\!\!\!\!\int A(x_t) A(y_t) d(\nu_2\times\nu_2)(x_{-n}, y_{-n})-
\\
2 \int\!\!\!\!\int A(x_t) A(y_t) d(\nu_1\times\nu_2)(x_{-n}, y_{-n})
\leq \Const \rho_2^n \enal \eneq since the first two terms are
exponentially close to $m_2(A\times A)$ and the last term is
exponentially close to $2 m_2(A\times A).$ Thus the result follows
by Borel--Cantelli. Q.E.D.

\begin{lemma}
$\lim_{n\to\infty} \mu_t^{(n)}(\nu)$ exists almost surely and if
$\nu_1$ and $\nu_2$ are two different measures then
$\lim_{n\to\infty} \mu_t^{(n)}(\nu_1)=
\lim_{n\to\infty} \mu_t^{(n)}(\nu_2)$ almost surely.
\end{lemma}
\begin{remark}
This proves part (A) of Theorem \ref{dissip}.
\end{remark}
\proof The second part follows immediately from Lemma \ref{LmBackMix}.
The first part follows the fact that almost surely there is a random
constant $C(\theta)$ such that
\begin{equation}
\label{Cauchy}
\left| \mu_t^{(n+1)}(A)-\mu_t^{(n)}(A)\right|\leq C(\theta) \rho^n
\end{equation}
The proof of (\ref{Cauchy}) is similar to the proof of (\ref{BackMix})
and can be left to the reader. Q.E.D.

By the construction we have part (B) of Theorem \ref{dissip}
\blm
With the notations above the family of measures
$\{\mu_t\}_{t\in\R_+}$ satisfies $\phi_{s,t} \mu_s=\mu_t$ almost
surely. \elm Thus, $t\mapsto \mu_t$ is a Markov process and it has
a continuous modification. Notice also that for any smooth
function $A$ on $M$ we have $\EXP(\mu_t(A))=m(A)$ for every $t\in
\R_+$, because $\EXP(\mu^{(n)}_{t}(A)) \to m(A)$ as $n \to
\infty$.

\blm With the notations above for any measure $\nu$ with
finite $p$-energy for some $p$ we have
$\forall A\in C^\infty(M)$
$$
\left|\nu(A(x_t))-\mu_t(A)\right|\leq C(\theta) e^{-\gamma t}.
$$
\elm
\begin{remark}
This is part (C) of Theorem \ref{dissip}.
\end{remark}
\proof Following the argument of the proof of Lemma \ref{LmBackMix}
we get for any two measures $\nu_1$ and $\nu_2$ of positive $p$-energy
$$
\PROB\left\{\left|\nu_1(A(x_t))-\nu_2(A(x_t))\right|\geq r\right\}
\leq \frac{\Const e^{-\gamma_1 t}}{r}.
$$
Taking $\nu_2=\mu_0$ we get
$$
\PROB\left\{\left|\nu(A(x_t))-\mu_t(A)\right|\geq r\right\}\leq
\frac{\Const e^{-\gamma_1 t}}{r}.
$$
The rest of the proof is similar to the proof of Lemma \ref{CrEquid}.
Q.E.D.

{\it Proof of Theorem \ref{ThDCLT}:}
The proof of this Theorem is the same as
the proof of Theorem \ref{measure}. Q.E.D.

\begin{remark}
In the conservative case $\mu_t\equiv m,$ so $\bar a$ and $\bralpha_k$
are
non-random and so Theorem \ref{ThDCLT} reduces to Theorem \ref{measure}.
Conversely, if $\forall A\in C^\infty(M)$ we have
$D''(A)=0$ then $\mu_t$ does not depend on $t$ and do
$\nu_t(A)=\EXP(\mu_t(A))=m(A).$ Using (\ref{BackVar}) we get
$$\EXP\left(\left[\mu_t^{(n)}(\nu)(A)-m(A)\right]^2\right)=$$
$$\int\!\!\!\!\int A(x_t) A(y_t) d(\nu\times\nu)(x_{-n}, y_{-n})+
m(A)^2-
2 m(A) \int A(x_t) d(x_{-n})\leq
\Const \rho_2^n$$
The first term here is close to $m_2(A\times A)$ and the sum of the other
two is close to $-m(A)^2.$ Since $n$ is arbitrary we get
$m_2(A\times A)=m(A)^2.$ By polarization $\forall A, B\in C(M)$
$$ \int\!\!\!\! \int A(x) B(y) dm_2(x,y)=m(A)m(B).$$
Hence $m_2=m\times m.$ By \cite{K} this implies that
each $\phi_{s,t} $ preserves $m.$ Hence we can characterize the
conservative case by the condition that the drift in
Corollary \ref{dis-measure} is non-random.
\end{remark}

{\bf Acknowledgement.}
This paper was started during the conference Nonlinear Analysis, 2000.
D. D. is grateful to IMA for the travel grant to attend this conference.
During the work on this paper we enjoyed the hospitability of IMPA,
Rio de Janeiro and of Universidad Autonoma de Madrid. We are espesially
grateful to Antonio Cordoba and Diego Cordoba for providing an excellent
working conditions in Madrid. Authors would like to thank P. Friz and
S.R.S. Varadhan for useful comments.


\begin{thebibliography}{9999}
\def\bitem#1{\bibitem[#1]{#1}}

\bitem{Ba} P.  Baxendale, {\it The Lyapunov spectrum of a stochastic
flow of diffeomorphisms}, LNM,  {\bf 1186}, 322--337, (1986);

\bitem{BH} P. Baxendale, T. Harris, {\it Isotropic Stochastic Flows},
Ann. Prob., {\bf 14}, 1155--1179, (1986);

\bitem{BS} P.  Baxendale, D. W. Stroock, {\it Large deviations and
stochastic flows of
diffeomorphisms,} Prob. Th. \& Rel. Fields {\bf 80}, 169--215, (1988);

\bitem{Bk} Gar. Birkhoff, {\it Lattice Theory} Providence, 1967;

\bitem{BL} P. Bougerol, J. Lacroix, {\it Products of random matrices
with applications to Schrodinger operators,} 1985;

\bitem{CC} R. Carmona, F. Cerou, {\it Transport by incompressible
random velocity fields: Simulations and Mathematical Conjectures},
153--181, AMS, Providence, 1999;

\bitem{C1} A. Carverhill, {\it A Formula for the Lyapunov numbers of
a stochastic flows. Applications to a perturbation theorem},
Stochastics, {\bf 14}, 209--225, (1985);

\bitem{C2} A. Carverhill, {\it Flows of stochatic dynamical systems:
ergodic theory}, Stochastics, {\bf 14}, 273--317, (1985);


\bitem{CSS1} M. Cranston, M. Scheutzow, D. Steinsaltz,
{\it Linear bounds for stochastic dispersion}, Ann. Probab.
{\bf 28}, no.4, 1852--1869, (2000);

\bitem{CSS2} M. Cranston, M. Scheutzow, D. Steinsaltz,
{\it Linear and near--linear bounds for stochastic dispersion,}
Elect. Comm. in Probability, {\bf 4}, 91--101, (1999);

\bitem{CS} M. Cranston, M. Scheutzow, {\it Dispersion
rates under finite mode Kolmogorov flows}, preprint;

\bitem{Da} R. Davis, {\it Lagrangian Ocean Studies},
Ann. Rev. Fluid ech., {\bf 23}, 43--64, (1991);

\bitem{Do} J. L. Doob, {\it Stochastic Processes}, John Wiley and
Sons;

\bitem{ES} K. Elworthy, D. Stroock, {\it Large deviation theory
for mean exponents of stochastic flows}, Stochastic Processes,
mathematics and physics, LNM, {\bf 1158}, 72--80, (1986);

\bitem{H} T. Harris, {\it Brownian Motion on the homeomorphisms of
the plane}, Ann. Prob., {\bf 9}, 232--254, (1981);


\bitem{I} M. Isichenko, {\it Percolation, Statistical topography
and transport in random media}, Reviews in Modern Physics, 1992;

\bitem{K} H. Kunita, {\it Stochastic Flows and Stochastic
Differential Equations}, Cambridge University Press, Cambridge, 1990;

\bitem{LY1} F. Ledrappier \& L.-S. Young,
{\it Dimension formula for random transformations,}
Comm. Math. Phys. {\bf 117}, 529--548, (1988);

\bitem{LY2}
F. Ledrappier \& L.-S. Young,
{\it Entropy formula for random transformations,}
Prob. Th., Rel. Fields {\bf 80}, 217--240, (1988);

\bitem{L1} Y. Le Jan, {\it On Isotropic Brownian Motion},
Zeit. fur Wahr., {\bf 70}, 609--620, (1985);

\bitem{L2} Y. Le Jan, {\it Hausdorff dimension for the statistical
equilibrium of stochastic flows,} Lect. Notes in Math. {\bf 1158},
201--207, (1986);

\bitem{L3} Y. Le Jan, {\it Equilibre statistique pour les produits
de diffeomorphisms aleatores independants,} C.R. Acad. Sci. Par.
Ser. I Math. {\bf 302}, (1986), 351--354;

\bitem{Le} R. Leandre, {\it Minoration en temps petit de la
densite d'une diffusion degeneree}, Journal of Functional
Analysis, {\bf 74}, 399-414, (1987);

\bitem{LS} H. Lisei, M. Scheutzow, {\it Linear Bounds and
Gaussian Tails in a Stochastic Dispersion Model}, preprint;

\bitem{Ma} P. Mattila, {\it Geometry of Sets and Measures in Euclidean
Spaces}, Cambridge University Press, 1995;

\bitem{Mo} S. Molchanov, {\it Topics in Statistical Oceanography:
in Stochastic Modelling in Physical Oceanography}, 343--380,
Boston, 1996, Birkhauser;

\bitem{MY} A. Monin, A. Yaglom, {\it Statistical Fluid Mechanics:
Mechanism of Turbulence}, MIT Press, Cambridge, MA, 1971;

\bitem{RT} D. Ruelle, F. Takens {\it On the nature of turbulence,}
Comm. Math. Phys. {\bf 20}, 167--192, (1971);

\bitem{SS} M. Scheutzow, D. Steinzaltz, {\it Chasing balls through
martingale fields}, preprint;

\bitem{Y} A. Yaglom, {\it Correlation Theory of Stationary and Related
Random Functions }, vol. 1: Basic Results, Springer-Verlag, New
York, 1987.

\bitem{ZC} C. Zirbel, E. Cinlar, {\it Dispertion of particles in
Brownian flows}, Adv. Appl. Prob., {\bf 28}, 53--74, (1996).

\end{thebibliography}
\end{document}